\documentclass[a4paper]{amsart}

\usepackage[utf8]{inputenc}
\usepackage[english]{babel}
\usepackage{bbm} %
\usepackage[tbtags]{mathtools}
\usepackage{amssymb,amsmath} %
\usepackage{color} %
\usepackage{graphicx} %
\usepackage{hyperref} %
\usepackage{cleveref}
\usepackage{booktabs} %
\usepackage{dsfont} %
\usepackage{tikz} %
\usepackage{a4wide}

\usepackage{pgfplots} 
\pgfplotsset{compat=newest}
\pgfplotsset{plot coordinates/math parser=false}
\newlength\figureheight
\newlength\figurewidth

\DeclareMathOperator{\Id}{Id}
\DeclareMathOperator{\diag}{diag}

\DeclareMathOperator{\appr}{appr}
\DeclareMathOperator{\MC}{MC}

\newcommand{\ul}[1]{\underline{#1}}

\newcommand{\dif}{ \ \mathrm{d}}

\newcommand{\dd}{\operatorname{d}\!}

\newcommand{\Cov}{\mathrm{Cov}}
\newcommand{\Cor}{\mathrm{Cor}}

\newcommand{\bfu}{\mathbf{u}}
\newcommand{\bfv}{\mathbf{v}}
\newcommand{\bfw}{\mathbf{w}}
\newcommand{\bfx}{\mathbf{x}}
\newcommand{\bfy}{\mathbf{y}}

\newcommand{\bfA}{\mathbf{A}}
\newcommand{\bfB}{\mathbf{B}}
\newcommand{\bfC}{\mathbf{C}}

\newcommand{\bfH}{\mathbf{H}}
\newcommand{\bfI}{\mathbf{I}}

\newcommand{\bfL}{\mathbf{L}}
\newcommand{\bfM}{\mathbf{M}}
\newcommand{\bfN}{\mathbf{N}}

\newcommand{\bfQ}{\mathbf{Q}}

\newcommand{\bfU}{\mathbf{U}}
\newcommand{\bfV}{\mathbf{V}}

\newcommand{\bbR}{\mathbb{R}}

\newcommand{\bflambda}{\boldsymbol{\lambda}}

\newcommand{\bfphi}{\boldsymbol{\phi}}

\newcommand{\bfzeta}{\boldsymbol{\zeta}}
\newcommand{\bfxi}{\boldsymbol{\xi}}

\newcommand{\calL}{\mathcal{L}}
\newcommand{\frakF}{\mathfrak{F}}

\DeclareMathOperator{\bcurl}{{\mathbf{curl}}}
\DeclareMathOperator{\bgrad}{\!{\boldsymbol{\nabla}}}

\DeclareMathOperator{\LandauO}{\mathcal{O}}
\DeclareMathOperator{\Landauo}{\textit{o}}

\date{\today}

\usepackage{amsthm}
\title{On uncertainty quantification of eigenvalues and eigenspaces with higher multiplicity}
\date{\today}
\author{Jürgen Dölz}
\author{David Ebert}
\address{Institute for Numerical Simulation, Friedrich-Hirzebruch-Allee 7,	53115 Bonn, Germany}
\email{\{doelz,ebert\}@ins.uni-bonn.de}
\thanks{This work was partially funded by the Deutsche Forschungsgemeinschaft (DFG, German Research Foundation) -- project number 501419255. The authors also acknowledge the support by the DFG under Germany's Excellence Strategy
	-- project number 390685813.}
\newtheorem{theorem}{Theorem}
\newtheorem{lemma}[theorem]{Lemma}
\newtheorem{corollary}[theorem]{Corollary}

\newtheorem{remark}[theorem]{Remark}

\numberwithin{theorem}{section}

\usepgfplotslibrary{external}
\begin{document} %

\begin{abstract}
	We consider generalized operator eigenvalue problems in variational form with random perturbations in the bilinear forms. 
	This setting is motivated by variational forms of partial differential equations with random input data.
	The considered eigenpairs can be of higher but finite multiplicity. 
	We investigate stochastic quantities of interest of the eigenpairs and discuss why, for multiplicity greater than 1, only the stochastic properties of the eigenspaces are meaningful, but not the ones of individual eigenpairs. 
	To that end, we characterize the Fréchet derivatives of the eigenpairs with respect to the perturbation and provide a new linear characterization for eigenpairs of higher multiplicity. As a side result, we prove local analyticity of the eigenspaces.
	Based on the Fréchet derivatives of the eigenpairs we discuss a meaningful Monte Carlo sampling strategy for multiple eigenvalues and develop an uncertainty quantification perturbation approach. 
	Numerical examples are presented to illustrate the theoretical results.
\end{abstract}

\maketitle

\section{Introduction}

\subsection{Motivation}
Eigenvalue problems have manifold applications in engineering, physics and medicine. Examples include particle accelerators \cite{ABB+2000}, photonic crystals \cite{dorfler_photonic_2011}, quantum theory \cite{CohenQuantum}, and structural engineering \cite{Tho2018}. To reduce development costs during design processes, the numerical simulation and solution of such eigenvalue problems has become increasingly more important.
However, since real-world circumstances typically deviate in an unforeseeable and seemingly random fashion from the ideal computational environment, the real-world solutions of the eigenvalue problem also deviate in a random fashion from the simulation results. To further reduce development costs in manufacturing processes, such uncertainties need to be quantified.

Given some random parameters $\mu\in X$ and $\varepsilon\in Y$ in Banach spaces $X$, $Y$, we consider parameter dependent eigenvalue problems
\begin{equation} \label{eq:genEVP}
	\begin{aligned}
	&\text{find }(u,\lambda)\in V\times\mathbb{R} \text{ with } u\neq 0 \text{ such that } 
	a(u,v;\mu) 
	= \lambda \ b(u,v;\varepsilon) \text{ for all }v\in V,
	\end{aligned}
\end{equation}
where $V\subset H$ are Hilbert spaces with dense and compact embedding, $a(\cdot,\cdot\,;\mu):V\times V \rightarrow \mathbb{R}$ is a continuous, $V$-elliptic and symmetric bilinear form, and $b(\cdot,\cdot\,;\varepsilon):H\times H \rightarrow \mathbb{R}$ is a continuous and $H$-elliptic scalar product on $H$. The dependence of the bilinear forms on $\mu$ and $\varepsilon$ is assumed to be linear, which is motivated by bilinear forms originating from the weak formulation of partial differential equations, some examples of which are given below. For each choice of parameters the standard theory applies \cite{Bof2010}, stating that there are countably many, positive eigenvalues whose reciprocals accumulate at zero.

\subsection{Examples} \label{sec:examples}
The standard examples for our considerations are
\begin{enumerate}
	\item The diffusion equation with homogeneous boundary data, where $V=H_0^1(D)\subset L^2(D)=H$ in some domain $D\subset\bbR^d$ and
	\[
	a(u,v;\mu)=\int_D\langle\mu\bgrad u,\bgrad v\rangle_{\bbR^d}\dd x,
	\qquad
	b(u,v;\varepsilon)=\int_D\varepsilon uv\dd x,
	\]
	with $\mu\in L^{\infty}(D;\bbR^{d\times d})=X$ whose values are symmetric matrices with positive eigenvalues which are uniformly bounded from above and below in $D$ and $\varepsilon\in L^{\infty}(D)=Y$ uniformly bounded from above and below in $D$.
	\item Maxwell's eigenvalue problem equipped with perfect electrical conductor boundary conditions, where we have $V=\bfH_0(\bcurl,D)\subset[L^2(D)]^3=H$ on a simply connected domain $D\subset\bbR^3$ with Lipschitz boundary. The bilinear forms read
	\[
	a(u,v;\mu)=\int_D\langle\mu\bcurl u,\bcurl v\rangle_{\bbR^3}\dd x,
	\qquad
	b(u,v;\varepsilon)=\int_D\langle\varepsilon u,v\rangle_{\bbR^3}\dd x,
	\]
	with $\mu,\varepsilon\in L^{\infty}(D;\bbR^{d\times d})=X=Y$ whose values are symmetric matrices with positive eigenvalues which are uniformly bounded from above and below in $D$.
	\item The shell eigenvalue problem on the domain $D = [-1,1]\times [0,2\pi]$, periodic in the second coordinate with boundary $\partial D = \{(x,y)\in D \colon x = \pm 1\}$ and assuming a time harmonic displacement field, where $V = [H^1_0(D)]^5 \subset [L^2(D)]^5=H$, has bilinear forms
\begin{align*}
	a(u,v;\mu) &= t a_m(u,v;\mu) + t a_s(u,v;\mu) + t^3 a_b(u,v;\mu) , \\
	b(u,v) &= t b_l(u,v) + t^3 b_r(u,v) ,
\end{align*}
for fixed $t>0$.
The bilinear form $a$ is comprised of the bilinear forms for membrane, shear and bending potential energies where $\mu \in L^\infty(D;\bbR)$ is Young's modulus and the bilinear form $b$ is comprised of the bilinear forms for displacement and rotation. 
For details on the individual parts and an in depth discussion we refer to \cite{dVHP2008,HL2019a}.
\end{enumerate}
We note that in all cases the parameters $\mu$ and $\varepsilon$ are assumed to be random but to satisfy the given conditions. Thus, the solutions to \cref{eq:genEVP} are also random.

\subsection{Related work}
While the mathematical theory on uncertainty quantification for partial differential and operator equations has been extensively developed in the last two decades, we only mention \cite{GHO2017,lord2014introduction,Smi2013,Soi2017} and the references therein, the uncertainty quantification of eigenvalue problems seems to have received relatively little attention in the mathematical community. Most of the early research on stochastic eigenvalue problems has been done in structural analysis and aerospace engineering, where the first works used a Monte Carlo approach in \cite{CT1969,SA1972b} and a perturbation approach in \cite{CT1969}. We refer to the review papers \cite{AF2007,BR1988,PV1977} and the references therein for a detailed overview. Nevertheless, it seems that only eigenvalues of single multiplicity and the uncertainty quantification of rather small systems, $3\times 3$-matrices in most of the references, were discussed. For larger systems, severe computational challenges occur such that these early approaches become prohibitively expensive.

To overcome these limitations, recent developments in the mathematical community include the analysis of a sparse grid approach in \cite{AndreevSchwab}, stochastic collocation methods in \cite{GSH2023,HL2019a,GACS2019}, and quasi-Monte Carlo methods in \cite{GGK+2019a,GS2022,GS2022a,Ngu2022}. Except for \cite{GACS2019}, \cite{GSH2023}, and \cite{HL2019a}, the first two of which considers a tracking technique to detect crossings of eigenvalues, all of these methods consider the case of eigenvalues of multiplicity one. The second and third reference, \cite{GSH2023} and \cite{HL2019a}, identify eigenspaces as the quantity of interest and derive analyticity and convergence rates for a stochastic collocation scheme for affine-parametric operators parametrized over $[-1,1]^{\mathbb{N}}$, as well as a spectral inverse and subspace iteration. As an alternative to sampling-based methods, stochastic Galerkin methods where discussed in \cite{BOS2019,ES2019,GG2007,HL2019b,Wil2010}, with main emphasis on the acceleration of the eigensolvers. However, crossings, bifurcations, and a rigorous error analysis of the numerical approximations compared to the analytic reality do not appear to have been the subject of investigation. Finally, we remark that parametrized eigenvalue problems have recently gained a lot of interest in the model order reduction community, see \cite{ABB+2022} and the references therein.

The reasons for the focus on eigenpairs with single multiplicity become clear when looking at the properties of the parameter-to-eigenpair map. For example, it is well known that eigenvalues and corresponding eigenspaces depend continuously on the parameters \cite{Kat1995}. However, this fact does not exclude that the eigenvalue trajectories exhibit crossings or bifurcations, which can occur on random occasions if uncertainties are present in the problem. The challenges concerning crossings and bifurcations of eigenvalues are more clearly recognizable after studying the seminal works due to Rellich. He found in a series of articles \cite{Rel1937,Rel1937a,Rel1939,Rel1940,Rel1941}, summarized in \cite{Rel1969}, that the eigenpairs to an eigenvalue of higher multiplicity are not necessarily Fr\'echet differentiable if the eigenvalue problem depends on more than one real parameter. This makes clear that, for stochastic dimensions larger than one, derivative based eigenvalue tracking algorithms need to be applied with caution. Another result due to Rellich is that if the problem indeed depends on a single real parameter, then there exist locally differentiable trajectories of the eigenpairs. Much later, in \cite{Dai1989,Nel1976}, a constructive method for finding the derivatives was presented. One might therefore ask why the perturbation approach of \cite{CT1969} for eigenvalues of single multiplicity cannot simply be extended to eigenvalues of higher multiplicity, as it was successfully done for quite a few settings for partial differential equations \cite{CNT2021,CS2013,Doe2020,EJ2020,HPS2013}. The challenge here is that the characterization of the \emph{first} derivative due to Rellich requires \emph{second} derivatives of the input data. As a consequence, the  Fr\'echet derivative is characterized as a \emph{non-linear} mapping rather than a \emph{linear} mapping, which makes its analysis and computation more involved than desirable.

Summarizing, most methods for uncertainty quantification of eigenvalue problems to date are restricted to the case where the eigenpair belongs to an eigenvalue of single multiplicity or a \emph{nondegenerate} eigenvalue. Perturbation approaches for eigenpairs to multiple eigenvalues, or \emph{degenerate} eigenvalues, are hindered by the lack of suitable characterization of the Fr\'echet derivative.

\subsection{Contributions}
The contributions of this article are threefold.
\begin{enumerate}
\item We review the available results on eigenpair derivatives from the literature and reformulate them in general Hilbert spaces. After generalizing a result by Sun \cite{Sun1988} for finite dimensions to our setting, we characterize the Fr\'echet derivative of eigenpairs to eigenvalues of finite multiplicity as a \emph{linear} mapping given through the solution operator of a saddle point problem. As side result of the generalization, we prove that the eigenspaces of finite multiplicity depend analytically on the perturbation.
\item Based on the results of \cite{Dai1989} and Rellich, we discuss why computing statistical quantities of interest of an eigenfunction to an eigenvalue with higher multiplicity is generally not meaningful. Based on this discussion, we develop a meaningful strategy to relate samples in sampling-based methods and discuss a perturbation approach for the uncertainty quantification of eigenvalues and eigen\emph{spaces} of higher multiplicity.
\item We discuss the efficient numerical solution of the arising tensor product equations from the perturbation approach and compare our new perturbation approach to the Monte Carlo simulation and vice versa.
\end{enumerate}

\subsection{Outline}
The article is organized as follows. First, in \cref{sec:deterministic}, we provide a detailed functional analytic setting of the considered parametric eigenvalue problems and characterize the eigenpair Fr\'echet derivatives for the deterministic case. \Cref{sec:stochastic} is dedicated to the discussion of stochastic eigenvalue problems and suitable uncertainty quantification approaches. The Galerkin discretization of the derived formulas is discussed in \cref{sec:implementation}, which also discusses the efficient solution of the arising tensor product equations. \Cref{sec:experiments} provides numerical examples, leading to our conclusions in \cref{sec:conclusion}.

\section{Deterministic derivatives of eigenvalue problems} \label{sec:deterministic}

\subsection{Problem setting}\label{sec:detproblemsetting}
On a Hilbert space $U$ we introduce the Banach space of continuous bilinear forms
\[
B(U)=\big\{a\colon U\times U\to\bbR; \|a\|_{\text{op}}<\infty\big\}
\]
equipped with the operator norm
\[
\|a\|_{\text{op}}=\sup \limits_{u,v\neq 0} 
\frac{a(u,v)}{\|u\|_U \|v\|_U}.
\]
Let $V$, $H$ be real Hilbert spaces equipped with inner products $\langle\cdot,\cdot\rangle_V$ and $\langle\cdot,\cdot\rangle_H$ and with dense and compact embedding $V\hookrightarrow H$ and $X$, $Y$ Banach spaces. We consider the parametrized, general eigenvalue problem
\begin{align}\label{eq:vareigprob}
\text{find}~(u,\lambda)\in V\times\mathbb{R} \text{ with } u\neq 0 \text{ such that } 
a(u,v;\mu) 
= \lambda\, b(u,v;\varepsilon) \text{ for all }v\in V,
\end{align}
where $\mu\in X$ and $\varepsilon\in Y$,
\[
a\in\calL(X;B(V)),\qquad b\in\calL(Y;B(H)),
\]
and $a(\cdot,\cdot\,;\mu):V\times V \rightarrow \mathbb{R}$ are continuous, $V$-elliptic, and symmetric bilinear forms and $b(\cdot,\cdot\,;\varepsilon):H\times H \rightarrow \mathbb{R}$ are continuous and $H$-elliptic scalar products on $H$. We note that $(u,\lambda)$ depends on $\mu$ and $\varepsilon$ and more general assumptions could be made, but we stick with the current setting to keep exposition simple and avoid cluttering of notation.

The spectral theorem for symmetric and compact operators yields that under these assumptions, for every choice $\mu \in X,\varepsilon\in Y$, \cref{eq:vareigprob} has countably many positive and real eigenvalues, which are unique when ordered in increasing order, cf.~\cite{Bof2010}. We note, however, that the eigenfunctions are only unique up to linear scaling and orthogonal transformations within the eigenspaces. To fix the scaling of the eigenfunctions up to the choice of sign we also require normalization of the eigenfunction by the scalar product, i.e. 
\begin{align}
b(u,u;\varepsilon)=1 . \label{eq:normalization}
\end{align}

\begin{remark}
We note that \cref{eq:vareigprob} is equivalent to the generalized eigenvalue problem
\begin{align}\label{eq:eig_operator}
\text{find}~(u,\lambda)\in V\times\mathbb{R} \text{ with } u\neq 0 \text{ such that } 
A_\mu u = \lambda\, B_\varepsilon u~\text{in}~V',
\end{align}
with parametrized, symmetric linear operators
\[
A_{\mu}\in\calL(V,V'),\qquad B_{\varepsilon}\in\calL(H,H),
\]

defined through
\begin{align*}
\langle A_\mu u,v\rangle_H&=a(u,v;\mu)\quad\text{for all}~u,v\in V,\\
\langle B_\varepsilon u,v\rangle_H&=b(u,v;\varepsilon)\quad\text{for all}~u,v\in H.
\end{align*}
We note that
$A_{\mu}$ is $V$-elliptic and $B_{\varepsilon}$ is $H$-elliptic and both operators depend depend linearly and continuously on $\mu\in X$ and $\varepsilon\in Y$. Thus, all of the following considerations for \cref{eq:vareigprob} also hold for \cref{eq:eig_operator}.
\end{remark}

\begin{remark}\label{rem:standard}
Using the notation from the previous remark, we readily remark that $A_\mu$ is self-adjoint with respect to the $H$-duality product and invertible. Thus, using compactness of the embedding $V\hookrightarrow H$, its inverse can be considered as a compact and symmetric operator $A_\mu^{-1}\colon H\to H$. This makes $A^{-1/2}$ well defined and, setting $v=A^{1/2}u$, implies that \cref{eq:vareigprob} and \cref{eq:eig_operator} are equivalent to
\begin{align}\label{eq:standard_operator_evp}
	\text{find}~(v,\lambda)\in H\times\mathbb{R} \text{ with } v\neq 0 \text{ such that } 
	T_{\mu,\varepsilon} v = \lambda^{-1}\, v~\text{in}~H,
\end{align}
with
\[
T_{\mu,\varepsilon}=A_\mu^{-1/2}B_\varepsilon A_\mu^{-1/2}\colon H\to H.
\]
Since $T_{\mu,\varepsilon}$ is positive definite, symmetric, and compact, the standard theory applies to \cref{eq:standard_operator_evp} and transfers also to \cref{eq:vareigprob} and \cref{eq:eig_operator}.
\end{remark}

For a Fr\'echet differentiable function $f\colon Z_1\to Z_2$ between two Banach spaces $Z_1$, $Z_2$, we denote its derivative at point $z_1\in Z_1$ in direction $h_1\in Z_1$ by $D_{z_1}^{h_1}f$. In the following, we assume that
\begin{subequations} \label{eq:perturbation}
\begin{align} 
\mu\colon E\supset B(0,\alpha_0)\to X,&\qquad\alpha\mapsto\mu_\alpha,\\
\varepsilon\colon F\supset B(0,\beta_0)\to Y,&\qquad\beta\mapsto\varepsilon_\beta,
\end{align}
\end{subequations}
for some $\alpha_0,\beta_0>0$ and Banach spaces $E$, $F$, and, for simplicity, that $\mu$ and $\varepsilon$ are continuously Fr\'echet differentiable, i.e.,
\[
\mu\in C^1\big(B(0,\alpha_0);X\big),\qquad\varepsilon\in C^1\big(B(0,\beta_0);Y\big). 
\]

Our aim in the following subchapters is then to characterize the derivatives of eigenpairs to arrive at series expansions like
\begin{subequations} \label{eq:series}
	\begin{align}
	\lambda
	&= \lambda_0 
	+ D_{(0,0)}^{(\alpha,\beta)}\lambda
	+ \Landauo \big(\|\alpha\|_E+\|\beta\|_F\big), \label{eq:lambda_series}	\\
	u
	&= u_0 
	+ D_{(0,0)}^{(\alpha,\beta)}u
	+ \Landauo \big(\|\alpha\|_E+\|\beta\|_F\big) ,
	\end{align}
\end{subequations}
for $\|\alpha\|_E,\|\beta\|_F\to 0$ with $(u_0,\lambda_0)$ being an eigenpair at the reference point $(\alpha,\beta)=(0,0)$. We note that $\Landauo \big(\|\alpha\|_E+\|\beta\|_F\big)$ can be replaced by $\LandauO \big((\|\alpha\|_E+\|\beta\|_F)^2\big)$ for $\mu\in C^2\big(B(0,\alpha_0);X\big)$, $\varepsilon\in C^2\big(B(0,\beta_0);Y\big)$ and $\alpha\in B(0,\alpha_0)$, $\beta\in B(0,\beta_0)$.

\subsection{Notation for eigenvalues and -spaces with higher multiplicity}
To deal with eigenvalues of (finite) multiplicity $m>1$, we may abuse the notation to read
\begin{align*}
	a(\bfu,\bfv;\mu)
	={}&
	\begin{bmatrix}
		a([\bfu]_1,[\bfv]_1;\mu)& \ldots & a([\bfu]_m,[\bfv]_1;\mu) \\
		\vdots & \ddots & \vdots \\
		a([\bfu]_1,[\bfv]_m;\mu)& \ldots & a([\bfu]_m,[\bfv]_m;\mu)
	\end{bmatrix}\in\bbR^{m\times m},
\end{align*}
for all row vectors $\bfu = [[\bfu]_1,\ldots,[\bfu]_m],\bfv = [[\bfv]_1,\ldots,[\bfv]_m] \in V^m$ and apply a similar notation for $b$.
We then consider the problem
\begin{equation} \label{eq:eig_bf}
	\text{find }({\bfu},\bflambda)\in V^m\times\bbR^{m\times m} \text{ such that }
	a({\bfu},\bfv;\mu)=b({\bfu},\bfv;\varepsilon)\cdot\bflambda
	\text{ for all }\bfv\in V^m
\end{equation}
with the normalization constraint
\begin{align}\label{eq:eig_bf_normalized}
b(\bfu,\bfu;\varepsilon)=\bfI.
\end{align}
It is clear that $\bflambda$ is a diagonal matrix if all elements of $\bfu$ are eigenfunctions, i.e., satisfy \cref{eq:vareigprob}. On the other hand, it is important to note that \cref{eq:eig_bf} can still hold, even if the single elements of $\bfu$ do not satisfy \cref{eq:vareigprob}. 
The following lemma states that if \cref{eq:eig_bf} and \cref{eq:eig_bf_normalized} hold, then each element of $\bfu$ is a linear combination of eigenfunctions and vice versa.
\begin{lemma}
\Cref{eq:eig_bf,eq:eig_bf_normalized} hold if and only if there is an orthogonal matrix $\bfQ\in\bbR^{m\times m}$ depending on $(\mu,\varepsilon)$ such that all $([\tilde{\bfu}]_i,\tilde{\lambda}_i)$, $i=1,\ldots,m$, $\diag(\tilde{\lambda}_1,\ldots,\tilde{\lambda}_m)=\bfQ^\intercal\bflambda\bfQ$, $\tilde{\bfu}=\bfu\bfQ$, satisfy \cref{eq:vareigprob,eq:normalization}.
\end{lemma}
\begin{proof}
We only need to show that \cref{eq:eig_bf} implies \cref{eq:vareigprob} after a suitable orthogonal transformation. To this end, we remark that setting $\bfv=\bfu$ in \cref{eq:eig_bf} and using \cref{eq:eig_bf_normalized} implies that $\bflambda=a(\bfu,\bfu;\mu)$ is symmetric and thus has a diagonal form $\diag(\tilde{\lambda}_1,\ldots,\tilde{\lambda}_m)=\bfQ^\intercal\bflambda\bfQ$. Now, since \cref{eq:eig_bf} is equivalent to 
\[
a(\bfu\bfQ,\bfv;\mu)=b(\bfu\bfQ,\bfv;\varepsilon)\cdot(\bfQ^\intercal\bflambda\bfQ)
\]
for all $\bfv\in V^m$, this yields the assertion.
\end{proof}

In the following, we will focus on the case where $\bflambda_0=\lambda_0\bfI\in\bbR^{m\times m}$ at $(\mu,\varepsilon)=(\mu_0,\varepsilon_0)$, i.e., $\lambda_0$ is an eigenvalue of multiplicity $m$ with $m$ $b(\cdot,\cdot\,;\varepsilon_0)$-orthonormal eigenfunctions $\bfu_0\in V^m$. Our aim is then to find vectorized versions of the series expansions \cref{eq:series}, i.e.
\begin{subequations}\label{eq:series_deg}
\begin{align}
	\bfu 
	&= \bfu_0 
	+ D_{(0,0)}^{(\alpha,\beta)}\bfu 
	+ \LandauO\big((\|\alpha\|_E+\|\beta\|_F)^2\big),\\
	\bflambda
	&= \bflambda_0 
	+ D_{(0,0)}^{(\alpha,\beta)} \bflambda
	+ \LandauO\big((\|\alpha\|_E+\|\beta\|_F)^2\big),\label{eq:series_deg_eig}
\end{align}
\end{subequations}
with $(\bfu,\bflambda)$ satisfying \cref{eq:eig_bf}. We note that in general varying $(\alpha,\beta)$ will split the eigenspace of the eigenvalue into eigenspaces of lower multiplicity and that without futher considerations $\bflambda$ will not necessarily stay diagonal while changing $\alpha$ and $\beta$. %
We will discuss the problem of finding a \emph{polarized} choice of the basis of the eigenspace to keep $\bflambda$ diagonal in \cref{sec:polarization} below.

\subsection{A regularity, analyticity, and orthogonality result}
Before characterizing the eigenpair derivatives in \cref{eq:series} and \cref{eq:series_deg} we first prove their existence. To this end, the proof is inspired by \cite{Sun1988} for the finite dimensional case of matrix eigenvalue problems. The following theorem is a generalization to the infinite dimensional case and proves analytic parametric regularity of the eigenvalues and eigenspaces as well as an orthogonality result which will become useful later on.

\begin{theorem}\label{thm:orthogonalitythm}
	Let $\lambda_0$ be an $m$-fold eigenvalue of \cref{eq:vareigprob} at $(\mu_0,\varepsilon_0)$ with eigenspace $U_0$ and $b(\cdot,\cdot\,;\varepsilon_0)$-orthonormal eigenbasis $\bfu_0$. Then there exists a unique local, analytic trajectory $(\mu,\varepsilon)\mapsto(\bfu,\bflambda)$ such that $(\bfu,\bflambda)$ satisfies \cref{eq:eig_bf} (with $\bflambda$ not necessarily being diagonal) and, at $(\mu_0,\varepsilon_0)$, it holds $(\bfu,\bflambda)=(\bfu_0,\lambda_0\bfI)$.
	
	Moreover, it holds $\bfu-\bfu_0\in U_0^\perp$, i.e., all elements of $\bfu-\bfu_0$ are locally $b(\cdot,\cdot\,;\varepsilon_0)$-orthogonal to $U_0$.
\end{theorem}
\begin{proof}
	For notational convenience during the proof we recall the one-to-one correspondence between the variational eigenvalue problem \cref{eq:vareigprob} and the generalized operator eigenvalue problem \cref{eq:eig_operator}.
	
	Let $U_{0}^\perp$ the $b(\cdot,\cdot\,;\varepsilon_0)$-orthogonal complement of $U_{0}$ in $V$. 
	We further split the operators into their actions on $U_{0}\oplus U_{0}^\perp=V$ onto $V'=U_{0}'\oplus(U_{0}^\perp)'$
	\begin{align*}
		A_\mu\colon U_{0}\oplus U_{0}^\perp\to (U_{0})'\oplus(U_{0}^\perp)',
		&\qquad
		(u_0,u_0^\perp)
		\mapsto
		\begin{bmatrix}
			A_\mu^{0,0} & A_\mu^{0,\perp} \\
			A_\mu^{\perp,0} & A_\mu^{\perp,\perp}
		\end{bmatrix}
		\begin{bmatrix}
			u_0\\
			u_0^\perp
		\end{bmatrix},\\
		B_\varepsilon\colon U_{0}\oplus U_{0}^\perp\to (U_{0}')\oplus(U_{0}^\perp)',
		&\qquad
		(u_0,u_0^\perp)
		\mapsto
		\begin{bmatrix}
			B_\varepsilon^{0,0} & B_\varepsilon^{0,\perp} \\
			B_\varepsilon^{\perp,0} & B_\varepsilon^{\perp,\perp}
		\end{bmatrix}
		\begin{bmatrix}
			u_0\\
			u_0^\perp
		\end{bmatrix},
	\end{align*}
	and readily remark that $A_\mu^{0,\perp}=(A_\mu^{\perp,0})'$ and $B_\varepsilon^{0,\perp}=(B_\varepsilon^{\perp,0})'$ and the diagonal blocks are self-adjoint in the $H$-inner product due to the self-adjointness of $A_\mu$ and $B_\varepsilon$.
	
	Similarly, we remark that a finite system of eigenvalues separated from the rest of the spectrum changes locally continuous under perturbation, see, e.g., \cite[Theorem 3.16 and Chapter 4.3.5]{Kat1995} applied to $T_{\mu,\varepsilon}$ from \cref{rem:standard}. I.e., there is a sufficiently small neighbourhood of $(\mu_0,\varepsilon_0)$ for which there is a continuous mapping $(\mu,\varepsilon)\mapsto(\bfu,\bflambda)$ with $(\mu_0,\varepsilon_0)\mapsto(\bfu_0,\lambda_0\bfI)$ and $(\bfu,\bflambda)$ being the solution to \cref{eq:eig_bf}. Moreover, in this neighbourhood of $(\mu_0,\varepsilon_0)$, there are, counting multiplicities, exactly $m$ solutions to \cref{eq:vareigprob} with an eigenvalue in a neighbourhood of $\lambda_0$.
	We denote the $m$-dimensional space spanned by the corresponding eigenfunctions as $U$ and its $b(\cdot,\cdot\,;\varepsilon)$-orthogonal complement as $U^\perp$. The $b(\cdot,\cdot\,;\varepsilon)$-orthogonality of $U$ and $U^\perp$ implies that $A_\mu$ and $B_\varepsilon$ have at least one block-diagonal representation in $U\oplus U^\perp$ and we claim that in a neighbourhood of $(\mu_0,\varepsilon_0)$ there exist
	\[
	W_{\mu,\varepsilon}\in\calL\big((U_0)',(U^\perp)'\big),
	\qquad
	Z_{\mu,\varepsilon}\in\calL\big(U,U_0^\perp\big),
	\]
	depending on $\mu$ and $\varepsilon$ such that
	\begin{subequations}\label{eq:ABblockdiagonal}
		\begin{align}
			\begin{bmatrix}
				(P^0)' & (Z_{\mu,\varepsilon})'\\
				W_{\mu,\varepsilon} & (P^\perp)'
			\end{bmatrix}
			\begin{bmatrix}
				A_\mu^{0,0} & A_\mu^{0,\perp} \\
				A_\mu^{\perp,0} & A_\mu^{\perp,\perp}
			\end{bmatrix}
			\begin{bmatrix}
				P^0 & (W_{\mu,\varepsilon})'\\
				Z_{\mu,\varepsilon} & P^\perp
			\end{bmatrix}
			&=
			\begin{bmatrix}
				A_{\mu,\varepsilon}^{(1)} & 0\\
				0 & A_{\mu,\varepsilon}^{(2)}
			\end{bmatrix}
			\colon
			U\oplus U^\perp\to (U)'\oplus (U^\perp)',\\
			\begin{bmatrix}
				(P^0)' & (Z_{\mu,\varepsilon})'\\
				W_{\mu,\varepsilon} & (P^\perp)'
			\end{bmatrix}
			\begin{bmatrix}
				B_\varepsilon^{0,0} & B_\varepsilon^{0,\perp} \\
				B_\varepsilon^{\perp,0} & B_\varepsilon^{\perp,\perp}
			\end{bmatrix}
			\begin{bmatrix}
				P^0 & (W_{\mu,\varepsilon})'\\
				Z_{\mu,\varepsilon} & P^\perp
			\end{bmatrix}
			&=
			\begin{bmatrix}
				B_{\mu,\varepsilon}^{(1)} & 0\\
				0 & B_{\mu,\varepsilon}^{(2)}
			\end{bmatrix}
			\colon
			U\oplus U^\perp\to (U)'\oplus (U^\perp)',
		\end{align}
	\end{subequations}
	with the $H$-orthogonal projections $P^0\colon U\to U_0$ and $P^\perp\colon U^\perp\to U_0^\perp$ onto $U_0$ and $U_0^\perp$ and
	\begin{align*}
		A_{\mu,\varepsilon}^{(1)}&=(P^0)'A_\mu^{0,0}P^0+(Z_{\mu,\varepsilon})'A_\mu^{\perp,0}P^0+(P^0)'A_\mu^{0,\perp}Z_{\mu,\varepsilon}+(Z_{\mu,\varepsilon})'A_\mu^{\perp,\perp}Z_{\mu,\varepsilon},\\
		A_{\mu,\varepsilon}^{(2)}&=W_{\mu,\varepsilon}A_\mu^{0,0}(W_{\mu,\varepsilon})'+W_{\mu,\varepsilon}A_\mu^{0,\perp}P^\perp+(P^\perp)'A_\mu^{\perp,0}(W_{\mu,\varepsilon})'+(P^\perp)'A_\mu^{\perp,\perp}P^\perp,\\
		B_{\mu,\varepsilon}^{(1)}&=(P^0)'B_\mu^{0,0}P^0+(Z_{\mu,\varepsilon})'B_\mu^{\perp,0}P^0+(P^0)'B_\mu^{0,\perp}Z_{\mu,\varepsilon}+(Z_{\mu,\varepsilon})'B_\mu^{\perp,\perp}Z_{\mu,\varepsilon},\\
		B_{\mu,\varepsilon}^{(2)}&=W_{\mu,\varepsilon}B_\mu^{0,0}(W_{\mu,\varepsilon})'+W_{\mu,\varepsilon}B_\mu^{0,\perp}P^\perp+(P^\perp)'B_\mu^{\perp,0}(W_{\mu,\varepsilon})'+(P^\perp)'B_\mu^{\perp,\perp}P^\perp.
	\end{align*}
	Exploiting that $A_\mu$ and $B_\varepsilon$ are self-adjoint we directly note that the necessary condition for $W_{\mu,\varepsilon}$ and $Z_{\mu,\varepsilon}$ for such a diagonal representation to hold is that the off-diagonal blocks of the matrix products must vanish. That is,
	\[
	\frakF(\mu,\varepsilon,W,Z)
	:=
	\begin{bmatrix}
		F(\mu,\varepsilon,W,Z)\\
		G(\mu,\varepsilon,W,Z)
	\end{bmatrix}
	=0
	\]
	where 
	\[
	F,G\colon X\times Y\times\calL\big((U_0)',(U^\perp)'\big)\times\calL\big(U,U_0^\perp\big)\to\calL\big(U,(U^\perp)'\big)
	\]
	are given by
	\begin{align*}
		F(\mu,\varepsilon,W,Z)&=(P^\perp)'A_\mu^{\perp,0}P^0+(P^\perp)'A_\mu^{\perp,\perp}Z+WA_\mu^{0,0}P^0+WA_\mu^{0,\perp}Z,\\
		G(\mu,\varepsilon,W,Z)&=(P^\perp)'B_\varepsilon^{\perp,0}P^0+(P^\perp)'B_\varepsilon^{\perp,\perp}Z+WB_\varepsilon^{0,0}P^0+WB_\varepsilon^{0,\perp}Z.
	\end{align*}
	Since $\frakF(\mu_0,\varepsilon_0,0,0)=0$ we have shown the claim if we can verify the assumptions of the Banach valued implicit function theorem, see, e.g., \cite[(10.2.1)]{Die1969}. To that end, we note that
	\[
	D_{(\mu_0,\varepsilon_0,0,0)}^{(0,0,\cdot,\cdot)}\frakF
	\colon
	\calL\big((U_0)',(U_0^\perp)'\big)\times\calL\big(U_0,U_0^\perp\big)
	\to
	\calL\big(U_0,(U_0^\perp)'\big)\times\calL\big(U_0,(U_0^\perp)'\big)
	\]
	with
	\[
	(\Xi,\Theta)
	\mapsto
	D_{(\mu_0,\varepsilon_0,0,0)}^{(0,0,\Xi,\Theta)}\frakF
	=
	\begin{bmatrix}
		\Xi A_{\mu_0}^{0,0} & A_{\mu_0}^{\perp,\perp}\Theta\\
		\Xi B_{\varepsilon_0}^{0,0} & B_{\varepsilon_0}^{\perp,\perp}\Theta
	\end{bmatrix}
	=
	\begin{bmatrix}
		\lambda_0\big[\cdot\,B_{\varepsilon_0}^{0,0}\big] & A_{\mu_0}^{\perp,\perp}\\
		\cdot\,B_{\varepsilon_0}^{0,0} & B_{\varepsilon_0}^{\perp,\perp}
	\end{bmatrix}
	\begin{bmatrix}
		\Xi\\
		\Theta
	\end{bmatrix},
	\]
	since $P^0$ acts as the identity on $U_0$ and $P^\perp$ as the identity on $U^\perp$. To show that
	$(\Xi,\Theta)\mapsto D_{(\mu_0,\varepsilon_0,0,0)}^{(0,0,\Xi,\Theta)}\frakF$
	is an isomorphism, we note that solving \cref{eq:vareigprob} at $(\mu_0,\varepsilon_0)$ is equivalent to computing the eigenpairs of the compact solution operator $S_{\mu_0,\varepsilon_0}=A_{\mu_0}^{-1}B_{\varepsilon_0}\colon H\to H$, which relate naturally to $T_{\mu_0,\varepsilon_0}$ from \cref{rem:standard}.
	Fredholm's alternative applied to $S_{\mu_0,\varepsilon_0}-1/\lambda_0$ implies that $A_{\mu_0}^{\perp,\perp}-\lambda_0B_{\varepsilon_0}^{\perp,\perp}$ is boundedly invertible which, using Gaussian elimination, shows that $(\Xi,\Theta)\mapsto D_{(\mu_0,\varepsilon_0,0,0)}^{(0,0,\Xi,\Theta)}\frakF$ is an isomorphism. Upon noting that $\mu\mapsto A_\mu$ and $\varepsilon\mapsto B_\varepsilon$ is linear, and thus analytic, the implicit function theorem yields that $W_{\mu,\varepsilon}$ and $Z_{\mu,\varepsilon}$ locally exist and are analytic.
	
	\Cref{eq:ABblockdiagonal} implies that every $u\in U$ can be represented as $u=u_0+Z_{\mu,\varepsilon}u$ for some $u_0\in U_0$. Thus, since $I-Z$ is invertible in a neighbourhood of $(\mu_0,\varepsilon_0)$, $\bfu=(I-Z_{\mu,\varepsilon})^{-1}\bfu_0$ is a basis of $U$ which is analytic in $(\mu,\varepsilon)$ for which it holds $\bfu-\bfu_0=Z_{\mu,\varepsilon}\bfu\in U_0^\perp$. Symmetry and ellipticity of $b(\cdot,\cdot\,;\varepsilon)$ imply that, in a neighbourhood of $(\mu_0,\varepsilon_0)$, $b(\bfu,\bfu;\varepsilon)$ is an invertible matrix such that analyticity of $\bflambda$ follows from testing \cref{eq:eig_bf} with $\bfu$ and solving for $\bflambda$.
\end{proof}
The proof of the theorem holds also for non-linear dependence of the bilinear forms on the parameters to obtain the following.
\begin{corollary}\label{cor:lowreg}
	Let the assumptions of \cref{thm:orthogonalitythm} hold, but with non-linear parameter dependence of the bilinear forms. Then, there hold the same implications as in \cref{thm:orthogonalitythm}, but the regularity of the mappings in \cref{thm:orthogonalitythm} and \cref{cor:eigenpair} is given by the combined regularity of the mapping $(\mu,\varepsilon)\mapsto\big(a(\cdot,\cdot\,;\mu),b(\cdot,\cdot\,;\varepsilon)\big)$.
\end{corollary}
A special case of this theorem was proven in \cite{GSH2023}, essentially extending the techniques from \cite{AndreevSchwab}. There, the parameter space was assumed to be $[-1,1]^{\mathbb{N}}$ and an analyticity result of the eigenspace was obtained by holomorphic extension. 

As a further special case we obtain the following result for nondegenerate eigenpairs, which, in various variants has already been proven in \cite{AndreevSchwab,Kat1995,Ngu2022,Rel1969}.
\begin{corollary}\label{cor:eigenpair}
	Let $(u_0,\lambda_0)$ be an eigenpair of \cref{eq:vareigprob} at $(\mu_0,\varepsilon_0)$ with nondegenerate eigenvalue. Then there exists a unique local, analytic trajectory $(\mu,\varepsilon)\mapsto(u,\lambda)$ such that $(u,\lambda)$ satisfies \cref{eq:vareigprob} and it holds $(u,\lambda)=(u_0,\lambda_0)$ at $(\mu_0,\varepsilon_0)$. If the parameter dependence of the bilinear forms is non-linear, the same implications as in \cref{cor:lowreg} hold.
\end{corollary}

\subsection{Derivatives of nondegenerate eigenvalues}
We start to characterize the derivatives of a single eigenvalue with multiplicity $m=1$. The case where $b$ is parameter independent can also be found in various textbooks on stability of matrix eigenvalue problems such as \cite{Saa2011} and is known as the Hellmann-Feynman theorem in quantum mechanics, cf. \cite[Chapter 5]{CohenQuantum}. While these techniques usually assume the existence of the involved derivatives, the existence is guaranteed by \cref{thm:orthogonalitythm} in our case.

\begin{lemma} \label{lem:lambda_derivatives}
Let $(u_0,\lambda_0)$ be a nondegenerate eigenpair of \cref{eq:vareigprob}. Then it holds
\begin{align}\label{eq:eigvalder}
D_{(0,0)}^{(\alpha,\beta)}\lambda
=
a\big(u_0,u_0;D_{0}^{\alpha}\mu\big)
-
\lambda_0
b\big(u_0,u_0;D_{0}^{\beta}\varepsilon\big).
\end{align}
\end{lemma}
\begin{proof}
We consider the derivative of the Rayleigh quotient
\begin{align*}
	D_{(0,0)}^{(\alpha,\beta)} \lambda 
	&= D_{(0,0)}^{(\alpha,\beta)}\frac{a(u,u;\mu)}{b(u,u;\varepsilon)}\\
	&= \frac{\Big(2 a\big(D_{(0,0)}^{(\alpha,\beta)} u,u_0;\mu_0\big)+a\big(u_0,u_0;D_{0}^{\alpha}\mu\big)\Big)b(u_0,u_0;\varepsilon_0)}
	{b(u_0,u_0;\varepsilon_0)^2}\\
	&\qquad\qquad -\frac{a(u_0,u_0;\mu_0)\Big(2 b\big(D_{(0,0)}^{(\alpha,\beta)} u,u_0;\varepsilon_0\big)+b\big(u_0,u_0;D_{0}^{\beta}\varepsilon\big)\Big)}
	{b(u_0,u_0;\varepsilon_0)^2}.
\end{align*}
The assertion follows from \cref{eq:vareigprob} and \cref{eq:normalization}. 
\end{proof}

\subsection{Derivatives of eigenfunctions to nondegenerate eigenvalues} \label{sec:eigenfunctions}
To determine the derivative of the eigenfunction we take the derivative of the equation in problem \cref{eq:vareigprob} to obtain
\begin{align}\label{eq:u1}
a\big(D_{(0,0)}^{(\alpha,\beta)} u,v;\mu_0\big)
- \lambda_0 b\big(D_{(0,0)}^{(\alpha,\beta)}u,v; \varepsilon_0\big)
=
\big(D_{(0,0)}^{(\alpha,\beta)}\lambda\big) b\big(u_0,v;\varepsilon_0\big)
- a\big(u_0,v;D_{0}^{\alpha}\mu\big)
+ \lambda_0 b\big(u_0,v;D_{0}^{\beta}\varepsilon\big).
\end{align}
Note that setting $v=u_0$ recovers \cref{eq:eigvalder}. Unfortunately, this characterization of $D_{(0,0)}^{(\alpha,\beta)} u{}\in V$ is not unique, since $D_{(0,0)}^{(\alpha,\beta)} u+cu_0$ solves \cref{eq:u1} for all $c\in\bbR$.
A unique representation of the eigenfunction derivative can be obtain by taking the derivative of the normalization condition \cref{eq:normalization}, which yields
\begin{align}\label{eq:u1cond}
0
=
2b\big(D_{(0,0)}^{(\alpha,\beta)} u,u_0;\varepsilon_0\big) + b\big(u_0,u_0;D_{0}^{\beta}\varepsilon\big).
\end{align}
Although we can consider $D_{(0,0)}^{(\alpha,\beta)}\lambda$ as known due to \cref{lem:lambda_derivatives}, solving for $(D_{(0,0)}^{(\alpha,\beta)}u,D_{(0,0)}^{(\alpha,\beta)}\lambda)$ simultaneously allows a simple characterization of the eigenpair derivative in terms of a saddle point problem.
To simplify the presentation, we introduce bilinear forms $A\colon V\times V\to\bbR$ and $B\colon\bbR\times V\to\bbR$ defined by
\begin{subequations}\label{eq:helpersSaddlepoint}
	\begin{align} 
		A(u,v)={}&a(u,v;\mu_0)-\lambda_0b(u,v;\varepsilon_0),\\
		B(\zeta,v)={}&\zeta b(u_0,v;\varepsilon_0).
	\end{align}
\end{subequations}
Solving \cref{eq:u1} with constraint \cref{eq:u1cond} is then equivalent to finding $(D_{(0,0)}^{(\alpha,\beta)}u,D_{(0,0)}^{(\alpha,\beta)}\lambda)\in V\times\bbR$ such that
\begin{align}\label{eq:saddlepoint}
\begin{aligned}
A\big(D_{(0,0)}^{(\alpha,\beta)}u,v\big)-B\big(D_{(0,0)}^{(\alpha,\beta)}\lambda,v\big)={}&-a\big(u_0,v;D_{0}^{\alpha}\mu)+\lambda_0 b\big(u_0,v;D_{0}^{\beta}\varepsilon\big),\\
B\big(\zeta,D_{(0,0)}^{(\alpha,\beta)}u\big)={}&-\frac{\zeta b\big(u_0,u_0;D_{0}^{\beta}\varepsilon\big)}{2},
\end{aligned}
\end{align}
for all $(v,\zeta)\in V\times\bbR$. We note that, in the finite dimensional case, the matrix formulation of these formulas coincides with \cite{Nel1976}, see also \cref{sec:implementation}.

\begin{lemma}\label{lem:solvenondegsaddle}
	The saddle point problem \cref{eq:saddlepoint} is uniquely solvable.
\end{lemma}
\begin{proof}
	
	We first show the LBB-condition of $B$. To this end, we estimate
	\[
	\inf_{0\neq\zeta\in\bbR}\sup_{0\neq v\in V}\frac{B(\zeta,v)}{|\zeta|\|v\|_V}
	=
	\inf_{0\neq\zeta\in\bbR}\sup_{0\neq v\in V}\frac{\zeta b(u_0,v;\varepsilon_0)}{|\zeta|\|v\|_V}
	\geq
	\frac{ b(u_0,u_0;\varepsilon_0)}{\|u_0\|_V}
	=
	\frac{1}{\|u_0\|_V}
	\geq
	C\frac{1}{\sqrt{\lambda_0}}>0,
	\]
	where we set $v=u_0$, used the normalization constraint $b(u_0,u_0;\varepsilon_0)=1$, and the $V$-ellipticity of $a(\cdot,\cdot;\mu_0)$ through
	\[
	\|u_0\|_V\leq C\sqrt{a(u_0,u_0;\mu _0)}=C\sqrt{\lambda_0b(u_0,u_0;\varepsilon_0)}=C\sqrt{\lambda_0}.
	\]
	
	It remains to show that finding $w\in W$ such that
	\begin{align}\label{eq:Ainvertibility}
	A(w,v)=\ell(v)
	\end{align}
	for all $v\in W$ with
	\[
	W=\ker B =\{v\in V\colon B(\zeta,v)=0~\text{for all}~\zeta\in\bbR\}=\operatorname{span}\{u_0\}^\perp
	\]
	is uniquely solvable for all $\ell\in W'$. Similar to the proof of \cref{thm:orthogonalitythm}, this follows from the Fredholm alternative applied to $S_{\mu_0,\varepsilon_0}-1/\lambda_0$ with $S_{\mu_0,\varepsilon_0}$ as in the proof of \cref{thm:orthogonalitythm}. Standard saddle point theory implies the assertion, see also \cite{BF1991}.
\end{proof}

\subsection{The problem with derivatives of degenerate eigenpairs} \label{sec:degen_eigenval}

Unfortunately, as we show in this subsection, following the same procedure to characterize derivatives of degenerate eigenpairs does not provide us with a satisfactory characterization of the derivatives. Taking the derivative of \cref{eq:eig_bf} yields
\begin{align}\label{eq:main_cond_deg}
\begin{aligned}
a\big(D_{(0,0)}^{(\alpha,\beta)} \bfu, \bfv;\mu_0\big)
- &\lambda_0 b\big(D_{(0,0)}^{(\alpha,\beta)}\bfu,\bfv; \varepsilon_0\big)\\
&\qquad=
b\big(\bfu_0,\bfv;\varepsilon_0\big)\big(D_{(0,0)}^{(\alpha,\beta)}\bflambda\big)
- a\big(\bfu_0,\bfv;D_{0}^{\alpha}\mu\big)
+ \lambda_0 b\big(\bfu_0,\bfv;D_{0}^{\beta}\varepsilon\big)
\end{aligned}
\end{align}
for all $\bfv \in V^m$, which coincides with \cref{eq:u1} for multiplicity $m=1$. Setting $\bfv=\bfu_0$ in \cref{eq:main_cond_deg}, this implies an analogous formula to \cref{eq:eigvalder}, i.e.,
\begin{align}\label{eq:degen_der_eigs}
D_{(0,0)}^{(\alpha,\beta)}\bflambda
=
a\big(\bfu_0,\bfu_0;D_{0}^{\alpha}\mu\big)
- \lambda_0 b\big(\bfu_0,\bfu_0;D_{0}^{\beta}\varepsilon\big).
\end{align}
The derivative of the orthonormality constraints 
of the eigenfunctions \cref{eq:eig_bf_normalized}
provides us with the normalization constraint
\begin{align}\label{eq:normcond_degen}
\mathbf{0}
=
b\big(D_{(0,0)}^{(\alpha,\beta)}\bfu,\bfu_0;\varepsilon_0\big)
+
b\big(\bfu_0,D_{(0,0)}^{(\alpha,\beta)}\bfu;\varepsilon_0\big)
+
b\big(\bfu_0,\bfu_0;D_{0}^{\beta}\varepsilon\big),
\end{align}
in complete analogy to equations \cref{eq:u1cond}.
We stress that $b(\bfu,\bfv;\varepsilon_0)\in\bbR^{m\times m}$ is not symmetric as a matrix for $m>1$, but that it holds
\begin{align*}
	b(\bfu,\bfv;\varepsilon_0) = b(\bfv,\bfu;\varepsilon_0)^\intercal
\end{align*}
due to the self-adjointness of the scalar product. Thus, the diagonal of \cref{eq:normcond_degen} implies
\begin{align}\label{eq:eigfuncderorth}
b\big(\big[D_{(0,0)}^{(\alpha,\beta)} \bfu\big]_i,[\bfu_0]_i;\varepsilon_0\big)=-\frac{b\big([\bfu_0]_i,[\bfu_0]_i;D_{0}^{\beta}\varepsilon\big)}{2},
\end{align}
for $i=1,\ldots,m$. The problem we face now is that the off-diagonal blocks of \cref{eq:normcond_degen} do not provide us with sufficiently many conditions to fully determine $D_{(0,0)}^{(\alpha,\beta)}\bfu$ for $m>1$. Thus, additional conditions need to be derived.

\subsection{Derivatives of eigenspaces to degenerated eigenvalues}
The traditional approach \cite{Dai1989,Rel1940} for obtaining additional constraints is to also consider second derivatives. In this approach, \cref{eq:eig_bf} is derived twice which, together with $\bfv=\bfu_0$ and \cref{eq:degen_der_eigs}, yields
\begin{align}\label{eq:constr_degen_nonlinear}
\begin{aligned}
&2a\big(D_{(0,0)}^{(\alpha,\beta)}\bfu,\bfu_0;D_{0}^{\alpha}\mu\big)
-
2\lambda_0 b\big(D_{(0,0)}^{(\alpha,\beta)}\bfu,\bfu_0;D_{0}^{\beta}\varepsilon\big)\\
&\qquad-
2\Big( b\big(D_{(0,0)}^{(\alpha,\beta)}\bfu,\bfu_0;\varepsilon_0\big)
+
b\big(\bfu_0,\bfu_0;D_{0}^{\beta}\varepsilon\big)\Big)D_{(0,0)}^{(\alpha,\beta)}\bflambda\\
&\qquad\qquad
=
\big(D_{(0,0)}^{(\alpha,\beta)}\big)^2\bflambda
-
a\big(\bfu_0,\bfu_0;\big(D_{0}^{\alpha}\big)^2\mu\big)
+
\lambda_0b\big(\bfu_0,\bfu_0;\big(D_{0}^{\beta}\big)^2\varepsilon\big),
\end{aligned}
\end{align}
where $\big(D_{(0,0)}^{(\alpha,\beta)}\big)^2$, $\big(D_{0}^{\alpha}\big)^2$, and $\big(D_{0}^{\beta}\big)^2$ denote the second Fr\'echet derivatives.
The idea is then that a suitable orthogonal transformation makes $\big(D_{(0,0)}^{(\alpha,\beta)}\big)^2\bflambda$ diagonal, and that the off-diagonal entries of \cref{eq:constr_degen_nonlinear} yield the additional constraints, see \cite{Dai1989} for details.

The problem for our purposes of uncertainty quantification is that this characterization of the Fr\'echet derivative $D_{(0,0)}^{(\alpha,\beta)}\bfu$ yields a \emph{non-linear} dependence on the perturbation terms $D_0^\alpha\mu$ and $D_0^\beta\varepsilon$. However, by the definition of Fr\'echet differentiability, there must be a \emph{linear} characterization and the second derivative should not be required. Thus, we will not follow the approach of using second derivatives. Instead, we use \cref{thm:orthogonalitythm} to derive a \emph{linear} dependence of the eigenspace derivatives on the perturbation parameters. To this end, \cref{thm:orthogonalitythm} directly implies that the derivative of an eigenfunction to an eigenvalue of higher multiplicity is orthogonal to the eigenspace while we neglect the normalization condition. Including the normalization condition leads to the following new characterization of the eigenspace derivatives which coincides with \cref{eq:saddlepoint} for $m=1$.

\begin{theorem} \label{th:Saddlepoint}
Let $\lambda_0$ be an eigenvalue of multiplicity $m$ at $(\mu_0,\varepsilon_0)$ of \cref{eq:vareigprob} with $b(\cdot,\cdot\,;\varepsilon_0)$-orthonormal eigenbasis $\bfu_0$. Let $(\mu,\varepsilon)\mapsto(\bfu,\bflambda)$ be the unique local, analytic trajectory such that $(\bfu,\bflambda)$ satisfies \cref{eq:eig_bf,eq:eig_bf_normalized} with coefficients \cref{eq:perturbation_linear} and, at $(\mu_0,\varepsilon_0)$, it holds $(\bfu,\bflambda)=(\bfu_0,\lambda_0\bfI)$. Let
\begin{subequations} \label{eq:helpersSaddlepointDeg}
\begin{align} 
A\colon V\times V&\to\bbR,&A(u,v)={}&a(u,v;\mu_0)-\lambda_0 b(u,v;\varepsilon_0),\\
B\colon \bbR^m\times V&\to\bbR,&B(\bfzeta,v)={}&\sum_{j=1}^m\zeta_j b([\bfu_0]_j,v;\varepsilon_0).
\end{align}
\end{subequations}
Then the derivatives $D_{(0,0)}^{(\alpha,\beta)}\bfu$ and $D_{(0,0)}^{(\alpha,\beta)}\bflambda$ are uniquely determined via the solutions of the saddle point problems ``Find $\big(\big[D_{(0,0)}^{(\alpha,\beta)}\bfu\big]_i,\big[D_{(0,0)}^{(\alpha,\beta)}\bflambda\big]_{:i}\big)\in V\times\bbR^m$ such that 
\begin{subequations} \label{eq:saddlepoint_bfu}
\begin{align}
	A\big(\big[D_{(0,0)}^{(\alpha,\beta)}\bfu\big]_i,v\big)-B\big(\big[D_{(0,0)}^{(\alpha,\beta)}\bflambda\big]_{:i},v\big)&=-a([\bfu_0]_i,v;D_{0}^{\alpha}\mu)+\lambda_0b([\bfu_0]_i,v;D_{0}^{\beta}\varepsilon),\\
	B\big(\bfzeta,\big[D_{(0,0)}^{(\alpha,\beta)} \bfu\big]_{i}\big)&=-\frac{\zeta_i b\big([\bfu_0]_i,[\bfu_0]_i;D_{0}^{\beta}\varepsilon\big)}{2},
\end{align}
\end{subequations}
holds for all $(v,\bfzeta)\in V\times\bbR^m$.'' Therein, by $\big[D_{(0,0)}^{(\alpha,\beta)}\bflambda\big]_{:i}$, we refer to the $i$.th column of $D_{(0,0)}^{(\alpha,\beta)}\bflambda\in\bbR^{m\times m}$.
\end{theorem}
\begin{proof}

The characterizing equations in \cref{eq:saddlepoint_bfu} follow from \cref{eq:main_cond_deg}, \cref{eq:eigfuncderorth}, and \cref{thm:orthogonalitythm}. It thus remains to show unique solvability of \cref{eq:saddlepoint_bfu}, which we prove by adapting the proof of \cref{lem:solvenondegsaddle} slightly. To this end, setting
\[
v_{\bfzeta}=\sum_{j=0}^m\zeta_j[\bfu_0]_j
\]
yields
\[
\|v_{\bfzeta}\|_V
\leq
C\sqrt{a(v_{\bfzeta},v_{\bfzeta};\mu_0)}
=
C\sqrt{\lambda_0b(v_{\bfzeta},v_{\bfzeta};\varepsilon_0)}
=
C\sqrt{\lambda_0\sum_{i,j=1}^m\zeta_i\zeta_jb([\bfu_0]_j,[\bfu_0]_i;\varepsilon)}
=
C\sqrt{\lambda_0}\|\bfzeta\|_{\bbR^m},
\]
with $C>0$ being the reciprocal ellipticity constant of $a(\cdot,\cdot;\mu_0)$. The LBB-condition of $B$ follows from
\begin{align*}
\inf_{0\neq\bfzeta\in\bbR^m}\sup_{0\neq v\in V}\frac{B(\bfzeta,v)}{\|\bfzeta\|_{\bbR^m}\|v\|_V}
&=
\inf_{0\neq\bfzeta\in\bbR^m}\sup_{0\neq v\in V}\frac{\sum_{j=1}^m\zeta_j b([\bfu_0]_j,v;\varepsilon_0)}{\|\bfzeta\|_{\bbR^m}\|v\|_V}\\
&\geq
\inf_{0\neq\bfzeta\in\bbR^m}\frac{\sum_{i,j=1}^m\zeta_j\zeta_i b([\bfu_0]_j,[\bfu_0]_i;\varepsilon_0)}{\|\bfzeta\|_{\bbR^m}\|v_{\bfzeta}\|_V}\\
&\geq
C\inf_{0\neq\bfzeta\in\bbR^m}\frac{\sum_{j=1}^m\zeta_j^2}{\sqrt{\lambda_0}\|\bfzeta\|_{\bbR^m}^2}\\
&\geq
C\frac{1}{\sqrt{\lambda_0}}>0.
\end{align*}
In analogy to the proof of \cref{lem:solvenondegsaddle} it remains to show that finding $w\in W$ such that $A(w,v)=\ell(v)$ holds for all $v\in W$ with
\[
W=\ker B =\{v\in V\colon B(\bfzeta,v)=0~\text{for all}~\bfzeta\in\bbR^m\}=\operatorname{span}\{[\bfu_0]_1,\ldots,[\bfu_0]_m\}^\perp
\]
is uniquely solvable for all $\ell\in W'$. This follows verbatim as in \cref{lem:solvenondegsaddle} from the Fredholm alternative.
\end{proof}

We note that, in contrast to \cref{thm:orthogonalitythm}, this new characterization of the eigenspace derivatives $D_{(0,0)}^{(\alpha,\beta)}\bfu$ is constructive, i.e., they are computable from any orthonormal basis $\bfu_0$ to $\lambda_0$.

\subsection{Polarization and derivatives of degenerated eigenvalues and eigenfunctions} \label{sec:polarization}

If the parameters $(\alpha,\beta)\in E\times F$ (and thus also $(\mu,\varepsilon)\in X\times Y$) in \cref{eq:vareigprob} depend on a single, real parameter $t$, it is known since \cite{Rel1937} that \cref{thm:orthogonalitythm} holds with some \emph{diagonal} ${\bflambda}$. More precisely, under the given assumptions, to an $m$-fold eigenvalue $\lambda_0$ with eigenspace $U_0$ there exist local, analytic trajectories $t\mapsto(\mu,\varepsilon)\mapsto([\tilde{\bfu}]_i,\tilde{\lambda}_i)$, $i=1,\ldots,m$, of orthonormal eigenpairs with $\tilde{\lambda}_i(0,0)=\lambda_0$. Given this existence result it remains to clarify how to find the corresponding eigenbasis of $U_0$ to these trajectories in a constructive way. The answer, up to second order, was given in \cite{Dai1989} whose basic idea is to exploit that \cref{eq:eig_bf} is invariant under orthogonal transformations of the eigenspace, i.e., for all solutions $(\bfu,\bflambda)$ and any orthogonal $\bfQ\in\bbR^{m\times m}$ it holds
\[
a({\bfu}\bfQ,\bfv;\mu)=b({\bfu}\bfQ,\bfv;\varepsilon)\cdot(\bfQ^\intercal\bflambda\bfQ)
\]
for all $\bfv\in V^m$.
Thus, assuming a fixed direction $(\alpha,\beta)\in E\times F$ for the moment, we may take any available $b(\cdot,\cdot\,;\varepsilon)$-orthonormal eigenbasis $\bfu$ of $U$ and look for an orthogonal $\bfQ\in\bbR^{m\times m}$, the \emph{polarization matrix}, such that $\bfu\bfQ$ coincides with the trajectories at $(\mu,\varepsilon)$. 
Substituting $\bfv=\bfu\bfQ$ in \cref{eq:main_cond_deg} and exploiting \cref{eq:degen_der_eigs} yields
\[
D_{(\alpha,\beta)}^{(\alpha,\beta)}\bflambda
=
a\big(\bfu,\bfu;D^{\alpha}_\alpha\mu\big)- b\big(\bfu,\bfu;D^{\beta}_\beta\varepsilon\big)\cdot\bflambda
=
\bfQ
\boldsymbol{\Lambda}
\bfQ^{\intercal}.
\]

Since $D_{(\alpha,\beta)}^{(\alpha,\beta)}\bflambda\in\bbR^{m\times m}$ is symmetric, it is clear that $D_{(\alpha,\beta)}^{(\alpha,\beta)}\bflambda$ can always be made diagonal by choosing $(\bfQ,\boldsymbol{\Lambda})$ as orthonormal eigenpairs of $D_{(\alpha,\beta)}^{(\alpha,\beta)}\bflambda$ with respect to the Euclidean inner product. 
The diagonal matrix $\bfQ^\intercal\bflambda\bfQ$ is unique up to permutation of the diagonal entries if the eigenvalues of $D_{(\alpha,\beta)}^{(\alpha,\beta)}\bflambda$ are distinct.

Assuming without loss of generality that the eigenvalues of $D_{(0,0)}^{(\alpha,\beta)}\bflambda$ are not degenerate themselves, we get the derivatives $D_{(0,0)}^{(\alpha,\beta)} \bfQ$ in analogy to \cref{eq:saddlepoint} by
\begin{align*}
	\big(D_{(0,0)}^{(\alpha,\beta)}\bflambda - [\boldsymbol{\Lambda}_0]_{ii} \bfI\big) [D_{(0,0)}^{(\alpha,\beta)}\bfQ]_i 
	&- [\bfQ_0]_i \frac{[D_{(0,0)}^{(\alpha,\beta)} \boldsymbol{\Lambda}]_{ii} }{2}	\\
	= 
	&-a\big((D_{(0,0)}^{(\alpha,\beta)}\bfu) [\bfQ_0]_i,\bfu_0;D_{0}^{\alpha} \mu \big) 
	 +b\big((D_{(0,0)}^{(\alpha,\beta)}\bfu) [\bfQ_0]_i,\bfu_0;D_{0}^{\beta} \varepsilon \big) \bflambda_0 \\
	&\quad-\frac{a\big(\bfu_0[\bfQ_0]_i,\bfu_0;(D_{0}^{\alpha})^2 \mu \big)}{2}
	 +\frac{b\big(\bfu_0[\bfQ_0]_i,\bfu_0;(D_{0}^{\beta})^2 \varepsilon \big) \bflambda_0}{2} \\
	&\quad+\frac{ b\big(\bfu_0[\bfQ_0]_i,\bfu_0;D_{0}^{\beta} \varepsilon \big) (D_{(0,0)}^{(\alpha,\beta)} \bflambda)}{2} \\
	[\bfQ_0]_i^\intercal [D_{(0,0)}^{(\alpha,\beta)} \bfQ]_i  
	=& \ 0,
\end{align*}
for $i=1,...,m$, where $\big[\big(D_{(0,0)}^{(\alpha,\beta)}\big)^2\bflambda\big]$ requires the existence of the second derivative of $\mu$ and $\varepsilon$ and $\boldsymbol{\Lambda}_0$ refers to the diagonal form of $D_{(0,0)}^{(\alpha,\beta)}\bflambda$.

Applying the found polarization matrix to \cref{eq:series_deg} yields
\begin{subequations}\label{eq:polarizedtaylor}
\begin{align}
\bfu\bfQ
&= \bfu_0 (\bfQ_0 + D_{(0,0)}^{(\alpha,\beta)} \bfQ)
+ (D_{(0,0)}^{(\alpha,\beta)}\bfu)\bfQ_0
+ \LandauO\big((\|\alpha\|_E+\|\beta\|_F)^2\big),\\
\bfQ^\intercal\bflambda\bfQ
&= \lambda_0\bfI
+ \bfQ_0^\intercal\big(D_{(0,0)}^{(\alpha,\beta)} \bflambda\big) \bfQ_0
+ \LandauO\big((\|\alpha\|_E+\|\beta\|_F)^2\big), 
\end{align}
\end{subequations}
where $(\bfu\bfQ,\bfQ^\intercal\bflambda\bfQ)=(\bfu\bfQ,\tilde{\boldsymbol{\lambda}})$ with $\tilde{\boldsymbol{\lambda}}$ diagonal solves \cref{eq:eig_bf}.

\begin{remark}
We emphasize that the above considerations hold for fixed $(\alpha,\beta)\in E\times F$ only, since the polarization matrix $\bfQ_0$ depends on $(\alpha,\beta)$ in general. Thus, a Taylor expansion can only be defined if the eigenproblem depends on a single, real parameter, see also \cite{Fri1996,Rel1937} for an in-depth discussion with examples and counter examples.
\end{remark}
\begin{remark}\label{rem:rewritefrechet}
In the special case where $E=F=\bbR$ we may rewrite the Fr\'echet derivative at $D_{(0,0)}^{(\alpha,\beta)}$ as
\begin{align*}
D_{(0,0)}^{(\alpha,\beta)} = \alpha D_{(0,0)}^{(1,0)} + \beta D_{(0,0)}^{(0,1)},
\end{align*}
such that the expansions read
\begin{align*}
\bflambda
&= \bflambda_0 
+ \alpha\big( D_{(0,0)}^{(1,0)}\bflambda\big) + \beta\big(D_{(0,0)}^{(0,1)} \bflambda\big)
+ \LandauO((|\alpha|+|\beta|)^2),\\
\bfu
&= \bfu_0 
+ \alpha\big( D_{(0,0)}^{(1,0)}\bfu\big) + \beta\big(D_{(0,0)}^{(0,1)}\bfu\big)
+ \LandauO((|\alpha|+|\beta|)^2).
\end{align*}
The polarization matrix $\bfQ_0$ becomes independent from $(\alpha,\beta)$ if $D_{(0,0)}^{(1,0)}\bflambda$ and $D_{(0,0)}^{(0,1)}\bflambda$ are diagonalizable over the same eigenbasis. The expansions then show that this is exactly the case if the eigenvalue trajectories can be approximated up to second order by planes in a neighbourhood of $\lambda_0$. Thus, in case of cone shaped trajectories as in the examples in \cref{sec:experiments}, a polarization matrix which is independent of $(\alpha,\beta)$ can never be found.
\end{remark}
We conclude that the eigenpairs to multiple eigenvalues have a G\^ateaux differentiable representation, but are not Fr\'echet differentiable in general.

\section{Uncertainty quantification of eigenpairs} \label{sec:stochastic}

\subsection{Problem setting}\label{sec:stochprob}
In the following, we aim to quantify the effect of uncertainties in the parameters $\mu$ and $\varepsilon$ on eigenpairs of finite multiplicity.
To that end, let $(\Omega,\mathcal{A},\mathbb{P})$ be a probability space and $\mu _0\in X$ and $\varepsilon_0\in Y$ deterministic reference parameters.
The parameters $\mu$ and $\varepsilon$ are then assumed to be modelled as random variables of the form
\begin{align} \label{eq:perturbation_linear}
	\mu(\omega,\alpha) 
	= \mu_0 + \alpha \mu_1(\omega) \in X, \qquad 
	\varepsilon(\omega,\beta) 
	= \varepsilon_0 + \beta \varepsilon_1(\omega) \in Y,
	\qquad\omega \in \Omega,
\end{align}
with $\alpha\in(-\alpha_0,\alpha_0)$, $\beta\in(-\beta_0,\beta_0)$ and $\mu_1 \in L_{\mathbb{P}}^2(\Omega;X)$, $\varepsilon_1 \in L_{\mathbb{P}}^2(\Omega;Y)$ uncorrelated and centered. We further assume $\|\mu_1\|_X \le C_\mu$ $\mathbb{P}$-a.s. and $\|\varepsilon_1\|_Y \le C_\varepsilon$ $\mathbb{P}$-a.s., for $C_\mu,C_\varepsilon <\infty$ fixed, which implies solvability of the eigenvalue problem \cref{eq:eig_bf} for all $\alpha,\beta$ if we choose $\alpha_0$ and $\beta_0$ sufficiently small. We note that this particular modelling corresponds to $E=F=\bbR$ in the setting of \cref{sec:deterministic}.

Rewriting the Fr\'echet derivative as in \cref{rem:rewritefrechet} yields
\[
D_{(0,0)}^{(\alpha,\beta)}\mu = D_0^\alpha \mu = \mu_1,
\qquad
D_{(0,0)}^{(\alpha,\beta)}\varepsilon = D_0^\beta \varepsilon = \varepsilon_1,
\]
and
\begin{subequations} \label{eq:series_stoch}
\begin{align}
	\bflambda(\omega)
	&= \bflambda_0 
	+ \alpha\big( D_{(0,0)}^{(1,0)}\bflambda\big) (\omega) + \beta\big(D_{(0,0)}^{(0,1)} \bflambda\big) (\omega)
	+ \LandauO((\alpha+\beta)^2,\omega)	
	\qquad \omega \in \Omega , \\
	\bfu(\omega) 
	&= \bfu_0 
	+ \alpha\big( D_{(0,0)}^{(1,0)}\bfu\big) (\omega) + \beta\big(D_{(0,0)}^{(0,1)}\bfu\big) (\omega)
	+ \LandauO((\alpha+\beta)^2,\omega)
	\qquad \omega \in \Omega,
\end{align}
\end{subequations}
due to \cref{eq:series_deg}. By $\LandauO(\cdot,\omega)$ we denote the usual Landau notation, but emphasize that the involved constant possibly depends on $\omega\in\Omega$.

In the following we aim to quantify the uncertainty in the eigenpairs $(\bfu(\omega),\bflambda(\omega))$. To that end, given two random fields $z_i\in L_{\mathbb{P}}^2(\Omega;Z_i)$, $i=1,2$, with $Z_1,Z_2$ being Hilbert spaces, we consider the mean and covariance
\begin{align*}
\mathbb{E}[z_1] 
= \int \limits_{\Omega} z_1 \dif \mathbb{P}\in Z_1, \qquad
\Cov[z_1,z_2]
= \mathbb{E}[(z_1-\mathbb{E}[z_1])\otimes(z_2-\mathbb{E}[z_2])]\in Z_1\otimes Z_2,
\end{align*}
with $\otimes$ denoting the Hilbertian tensor product, and abbreviate $\Cov[z_1]=\Cov[z_1,z_1]$. We start our discussion with a few fundamental considerations for eigenpairs of degenerate eigenvalues, shortly comment on sampling-based approaches, and then use the expansions in \cref{eq:series_stoch} to derive a perturbation approach.

\subsection{A perspective on the stochastic behaviour of degenerate eigenpairs} \label{sec:behav_stoch_deg}
For a meaningful uncertainty quantification of degenerate eigenpairs the different eigenpair realizations $\big(\bfu(\omega),\bflambda(\omega)\big)$ need to be related in a suitable sense for $\mathbb{P}$-a.e.\ $\omega\in\Omega$, in particular if crossings and bifurcations of the eigenvalue trajectories occur. In the following we use the perturbation formulas from the previous section as a theoretical tool to give a perspective on this subject. We emphasize that our findings also apply to other approaches, e.g.\ sampling-based approaches, for uncertainty quantification. To that end, pathwise application of our perturbation formulas \cref{eq:polarizedtaylor} and \cref{eq:series_stoch} yields

\begin{align*}
\bfQ(\omega)^\intercal\bflambda(\omega)\bfQ(\omega)
&= \bflambda_0 
+ \alpha\bfQ_0(\omega)^\intercal \big( D_{(0,0)}^{(1,0)}\bflambda\big) (\omega)\bfQ_0(\omega) 
+ \beta\bfQ_0(\omega)^\intercal \big(D_{(0,0)}^{(0,1)} \bflambda\big) (\omega)\bfQ_0(\omega)\\
&\hspace*{9cm}+ \LandauO\big((\alpha+\beta)^2,\omega\big),\\
\bfu(\omega)\bfQ(\omega)
&= \bfu_0 \Big( \bfQ_0(\omega) + \alpha\big(D_{(0,0)}^{(1,0)} \bfQ(\omega)\big) + \beta\big(D_{(0,0)}^{(0,1)} \bfQ(\omega)\big) \Big) \\
&\hspace*{3cm}
+ \alpha\big( D_{(0,0)}^{(1,0)}\bfu\big) (\omega)\bfQ_0(\omega) 
+ \beta\big(D_{(0,0)}^{(0,1)}\bfu\big) (\omega)\bfQ_0(\omega) \\
&\hspace*{9cm}
+ \LandauO\big((\alpha+\beta)^2,\omega\big),
\end{align*}
for all $\omega\in\Omega$. This relates all eigenpair realizations to an unperturbed, i.e., $\alpha=\beta=0$, reference eigenpair $(\bfu_0,\bflambda_0)$ with a fixed choice of basis and $\bflambda_0=\lambda_0\bfI$.

These expansions highlight the difficulties with relating different eigenpair realizations, most prominently visible by randomness of the polarization matrix $\bfQ_0(\omega)$. It causes the reference point for the eigenfunction expansion to be random and, since $\bfQ_0(\omega)$ is not unique for fixed $\omega$, causes the order of the eigenpairs to be random. Unfortunately, bare further assumptions, the randomness of the polarization matrix can not be avoided and any definition of stochastic quantities of interest must thus respect this randomness.

To that end, the seemingly only reasonable approach to identify different eigenpair realizations is to fix a \emph{reference polarization} into which all eigenpair realizations are transformed. This effectively means that we now need to perform uncertainty quantification for \emph{eigenspaces}. The following subsection elaborates on this viewpoint for sampling-based approaches. Afterwards we will use our expansions to develop a perturbation based approach.

\subsection{On sampling-based methods}\label{sec:sampling}
An intuitive way to compute an approximation of the expected value and the covariance of the $m$-fold eigenpair $(\bfu,\bflambda)$ are sampling-based approaches which are most prominently represented by the Monte Carlo approach. To that end, we draw realizations of the parameters $\mu_1(\omega)$ and $\varepsilon_1(\omega)$ to compute the eigenpairs $(\bfu(\omega),\bflambda(\omega))$.

Following the approach form \cref{sec:behav_stoch_deg} we pick the identity as a reference polarization such that we aim at relating $(\bfu(\omega),\bflambda(\omega))$ to $(\bfu_0,\bflambda_0)$, $\bflambda_0=\lambda_0\bfI$. Computing the singular value decomposition  
\[
b(\bfu(\omega),\bfu_0;\varepsilon_0) = \bfU(\omega) \boldsymbol\Sigma(\omega) \bfV(\omega)\in\bbR^{m\times m}
\]
implies $\boldsymbol\Sigma(\omega)\approx\bfI$ for $(\mu_1(\omega),\varepsilon_1(\omega))$ close to zero,
\[
b(\bfu(\omega)\bfV(\omega)^\intercal,\bfu_0\bfU(\omega)^\intercal;\varepsilon_0)
=
\bfU(\omega)^\intercal b(\bfu(\omega),\bfu_0;\varepsilon_0)\bfV(\omega)^\intercal
=
\boldsymbol\Sigma(\omega),
\]
and thus
\[
b(\bfu(\omega)\bfV(\omega)^\intercal\bfU(\omega),\bfu_0;\varepsilon_0)
=
b(\bfu(\omega),\bfu_0;\varepsilon_0)\bfV(\omega)^\intercal\bfU(\omega)
\approx
\bfI.
\]
This suggests that
\[
\big(\bfu(\omega)\bfV(\omega)^\intercal\bfU(\omega),\bfU(\omega)^\intercal\bfV(\omega)\bflambda(\omega)\bfV(\omega)^\intercal\bfU(\omega)\big)
\]
can be related to $(\bfu_0,\bflambda_0)$ without further consideration. We note that in many cases the computational cost of computing the SVD of an $m\times m$-matrix is negligible when compared to the original eigenvalue problem.
\begin{remark}
	Our approach can be understood as rotating the perturbed eigenspace onto the reference eigenspace. Another possibility, originally proposed in \cite{GSH2023,HL2019a}, would be to identify the eigenspaces by means of a spectral projection.
\end{remark}

\subsection{Series approximations of statistical quantities of interest}
Sampling approaches require the solution of an eigenvalue problem for every sample, which can be computationally costly or even prohibitive.
As an alternative, we derive in the following series expansions of the mean and the covariance of the eigenpairs.

The following two lemmata follow immediately from the assumptions in \cref{sec:stochprob}.
\begin{lemma} \label{lem:integrability}
	It holds  
	\begin{align*}
	\Big[\omega\mapsto a(\cdot,\cdot\,;\mu(\omega))\Big]
	\in L^2_\mathbb{P}(\Omega;B(V)), \qquad 
	\Big[\omega\mapsto b(\cdot,\cdot\,;\varepsilon(\omega))\Big]
	\in L^2_\mathbb{P}(\Omega;B(H)) .
	\end{align*}
\end{lemma}
\begin{proof}
	The definition of the Bochner spaces, $a\in\calL(X;B(V))$, and $\mu_1 \in L_{\mathbb{P}}^2(\Omega;X))$ imply
	\begin{align*}
	\int \limits_{\Omega} \| a(\cdot,\cdot\,;\mu(\omega)) \|^2_{B(V)} \dif \mathbb{P}(\omega) 
	\le \| a(\cdot,\cdot\,;\cdot) \|^2_{\mathcal{L}(X;B(V))}  
	\|\mu\|_{L_{\mathbb{P}}^2(\Omega;X)}^2<\infty.
	\end{align*} 
	The proof for $b$ is analogous.
\end{proof}
\begin{lemma} \label{lem:exp_bil}
	It holds $\mathbb{E}[a(u_0,v;\mu_1)] = 0$ and $\mathbb{E}[b(u_0,v;\varepsilon_1)] = 0$ for all $v\in V$.
\end{lemma}
\begin{proof}
	$a\in\calL(X;B(V))$ implies
	\begin{align*}
	\mathbb{E}[a(u_0,v;\mu_1(\omega))]
	= a(u_0,v;\mathbb{E}[\mu_1(\omega)]) 
	= a(u_0,v;0) 
	= 0\qquad\text{for all}~v\in V.
	\end{align*}
	The proof for b reads verbatim.
\end{proof}

With these preparatory lemmata in place we can state the main theorem of this section.
\begin{theorem} \label{th:stochastic}
	Let $\lambda_0$ be an eigenvalue of multiplicity $m$ at $(\mu_0,\varepsilon_0)$ of \cref{eq:vareigprob} with $b(\cdot,\cdot\,;\varepsilon_0)$-orthonormal eigenbasis $\bfu_0$. Let $(\mu,\varepsilon)\mapsto(\bfu,\bflambda)$ be the unique local, analytic trajectory such that $(\bfu,\bflambda)$ satisfies \cref{eq:eig_bf,eq:eig_bf_normalized} with coefficients \cref{eq:perturbation_linear} and, at $(\mu_0,\varepsilon_0)$, it holds $(\bfu,\bflambda)=(\bfu_0,\lambda_0\bfI)$. Then it holds
	\begin{align}
		\mathbb{E}[\bflambda] 
		&= \bflambda_0 
		+ \LandauO\big((|\alpha|+|\beta|)^2\big) , \label{eq:exp_lambda}\\
		\mathbb{E}[\bfu]
		&= \bfu_0
		+ \LandauO\big((|\alpha|+|\beta|)^2\big) , \label{eq:exp_u}\\		
		\Cov[\bflambda] 
		&= \alpha^2 \Cov[D_{(0,0)}^{(1,0)} \bflambda] 
		+ \beta^2 \Cov[D_{(0,0)}^{(0,1)} \bflambda]
		+ \LandauO\big((|\alpha|+|\beta|)^3\big) , \label{eq:cov_lambda} \\
		\Cov[\bfu]
		&= \alpha^2 \Cov[D_{(0,0)}^{(1,0)} \bfu] 
		+ \beta^2 \Cov[D_{(0,0)}^{(0,1)} \bfu] 
		+ \LandauO\big((|\alpha|+|\beta|)^3\big) . \label{eq:cov_u}
	\end{align}
\end{theorem}
\begin{proof}
We first remark that the assumptions of \cref{sec:stochprob} allow us to bound the $\LandauO((|\alpha|+|\beta|)^2,\omega)$-terms in \cref{eq:series_stoch} independently of $\omega$ $\mathbb{P}$-a.e. To that end, we remark that, for all $\omega\in\Omega$, $\alpha\mapsto\mu(\omega,\alpha)$ and $\beta\mapsto\varepsilon(\omega,\beta)$ are smooth functions which together with \cref{thm:orthogonalitythm} yields that $(\alpha,\beta)\mapsto(\bfu,\bflambda)$ is also smooth. Thus, the second derivative exists locally and its arguments are bounded due to $\|\mu_1\|_X \le C_\mu$ and $\|\varepsilon_1\|_Y \le C_\varepsilon$ for $\mathbb{P}$-a.e. Lagrange's form of remainder shows the claim.

The equations \cref{eq:exp_lambda} and \cref{eq:exp_u} for the mean follow by taking the mean of \cref{eq:series_stoch}, which yields
\begin{subequations}\label{eq:meanexpansion}
\begin{align}
	\mathbb{E}[\bflambda] 
	&= \bflambda_0 
	+ \alpha\mathbb{E}\big[D_{(0,0)}^{(1,0)}\bflambda\big] + \beta\mathbb{E}\big[D_{(0,0)}^{(0,1)}\bflambda\big]
	+ \LandauO\big((|\alpha|+|\beta|)^2\big) ,\\
	\mathbb{E}[\bfu] 
	&= \bfu_0 
	+ \alpha\mathbb{E}\big[D_{(0,0)}^{(1,0)}\bfu\big] + \beta\mathbb{E}\big[D_{(0,0)}^{(0,1)}\bfu\big]
	+ \LandauO\big((|\alpha|+|\beta|)^2\big).
\end{align}
\end{subequations}

Taking the mean of \cref{eq:saddlepoint_bfu} yields
\begin{align*}
A\big(\mathbb{E}\big[D_{(0,0)}^{(\alpha,\beta)}\bfu\big]_i,v\big)-B\big(\mathbb{E}\big[D_{(0,0)}^{(\alpha,\beta)}\bflambda\big]_{:i},v\big)&=0,\\
B\big(\bfzeta,\mathbb{E}\big[D_{(0,0)}^{(\alpha,\beta)} \bfu\big]_{i}\big)&=0,
\end{align*}
for all $(v,\bfzeta)\in V\times\bbR^m$, $i=1,\ldots,m$. Since the saddle point problem is uniquely solvable due to \cref{th:Saddlepoint}, this implies that the first derivatives in \cref{eq:meanexpansion} vanish.

For the covariance approximations \cref{eq:cov_lambda} and \cref{eq:cov_u}, we insert \cref{eq:series_stoch}, \cref{eq:exp_lambda}, and \cref{eq:exp_u} into the definition of the covariance to obtain
\begin{align*}
\Cov[\bflambda]
	&= \mathbb{E}[ \big(
	\bflambda - \mathbb{E}[\bflambda]
	\big) \otimes \big(
	\bflambda - \mathbb{E}[\bflambda]
	\big) ] \\
	&= \mathbb{E}
	[ \big(
	\alpha D_{(0,0)}^{(1,0)} \bflambda
	+ \beta D_{(0,0)}^{(0,1)} \bflambda 
	+ \LandauO\big((|\alpha|+|\beta|)^2\big) 
	\big) \otimes \big(
	\alpha D_{(0,0)}^{(1,0)} \bflambda
	+ \beta D_{(0,0)}^{(0,1)} \bflambda 
	+ \LandauO\big((|\alpha|+|\beta|)^2\big) 
	\big) ] \\
	&= \alpha^2 \Cov[D_{(0,0)}^{(1,0)} \bflambda] 
	+ \beta^2 \Cov[D_{(0,0)}^{(0,1)} \bflambda] 
	+ \LandauO\big((|\alpha|+|\beta|)^3\big) \\
	\Cov[\bfu] 
	&= \mathbb{E}[\big(
	\bfu-\mathbb{E}[\bfu]
	\big)\otimes\big(
	\bfu-\mathbb{E}[\bfu]
	\big)] \\
	&= \mathbb{E}[\big(
	\alpha D_{(0,0)}^{(1,0)} \bfu 
	+ \beta D_{(0,0)}^{(0,1)} \bfu	
	+ \LandauO\big((|\alpha|+|\beta|)^2\big)
	\big)\otimes \big(
	\alpha D_{(0,0)}^{(1,0)} \bfu 
	+ \beta D_{(0,0)}^{(0,1)} \bfu	
	+ \LandauO\big((|\alpha|+|\beta|)^2\big)
	\big)] \\
	&= \alpha^2 \Cov[D_{(0,0)}^{(1,0)} \bfu] 
	+ \beta^2 \Cov[D_{(0,0)}^{(0,1)} \bfu]
	+ \LandauO\big((|\alpha|+|\beta|)^3\big)
\end{align*}
For the last steps we note that the uncorrelatedness of $\mu_1$ and $\varepsilon_1$ implies
\[
\Cov\Big[D_{(0,0)}^{(1,0)}\bflambda,D_{(0,0)}^{(0,1)}\bflambda\Big]=0,
\qquad
\Cov\Big[D_{(0,0)}^{(1,0)}\bfu,D_{(0,0)}^{(0,1)}\bfu\Big]=0,
\]
and thus the assertion.
\end{proof}
We note that covariances between eigenvalues and eigenspaces as well as covariances across eigenspaces could be considered in complete analogy.
\begin{remark}
There are examples in the literature, such as partial differential equations on random domains modelled by the domain mapping approach \cite{HPS2016}, where $\mu$ and $\varepsilon$ are polynomials of higher degree in $\alpha$ and $\beta$ and correlated. In this case \cref{th:stochastic} and the following considerations hold with straightforward adaptions.
\end{remark}

\subsection{Covariance equations} \label{sec:variance}
For a numerical implementation of the covariance approximations \cref{eq:cov_lambda} and \cref{eq:cov_u} we need to characterize the arising correction terms. Adapting the approach from \cite{SchwabTodor} to saddle point problems we obtain the following.
\begin{theorem} \label{th:TensorVarDer}
	Let $\lambda_0$ be an eigenvalue of multiplicity $m$ at $(\mu_0,\varepsilon_0)$ of \cref{eq:vareigprob} with $b(\cdot,\cdot\,;\varepsilon_0)$-orthonormal eigenbasis $\bfu_0$. Let $(\mu,\varepsilon)\mapsto(\bfu,\bflambda)$ be the unique local, analytic trajectory such that $(\bfu,\bflambda)$ satisfies \cref{eq:eig_bf,eq:eig_bf_normalized} with coefficients \cref{eq:perturbation_linear} and, at $(\mu_0,\varepsilon_0)$, it holds $(\bfu,\bflambda)=(\bfu_0,\lambda_0\bfI)$. Then it holds
\begin{align*}
&\begin{bmatrix} 
A(\cdot,\bfv)\otimes\Id & -B(\cdot,\bfv)\otimes\Id \\
B(\bfzeta,\cdot)\otimes\Id & 0
\end{bmatrix}
\Cov
\begin{bmatrix}
D_{(0,0)}^{(1,0)} \bfu  \\
D_{(0,0)}^{(1,0)} \bflambda
\end{bmatrix}
\begin{bmatrix}
\Id\otimes A(\cdot,\bfw) & \Id\otimes B(\bfxi,\cdot) \\
-\Id\otimes B(\cdot,\bfw) & 0
\end{bmatrix} \\ & \hspace*{6.5cm}
= 
\Cov \begin{bmatrix} 
-a(\bfu_0,\cdot\,;\mu_1) \\ 0
\end{bmatrix}\big(\bfv\otimes\bfzeta,\bfw\otimes\bfxi\big), \\
&\begin{bmatrix} 
A(\cdot,\bfv)\otimes\Id & -B(\cdot,\bfv)\otimes\Id \\
B(\bfzeta,\cdot)\otimes\Id & 0
\end{bmatrix} \Cov
\begin{bmatrix}
D_{(0,0)}^{(0,1)} \bfu  \\
D_{(0,0)}^{(0,1)} \bflambda
\end{bmatrix}
\begin{bmatrix}
\Id\otimes A(\cdot,\bfw) & \Id\otimes B(\bfxi,\cdot) \\
-\Id\otimes B(\cdot,\bfw) & 0
\end{bmatrix} \\ & \hspace*{6.5cm}
= 
\Cov \begin{bmatrix} 
\lambda_0 b(\bfu_0,\cdot\,;\varepsilon_1) \\ -\diag \limits_{i=1,\ldots,m}\frac{\cdot \,b([\bfu_0]_i,[\bfu_0]_i;\varepsilon_1)}{2}
\end{bmatrix}\big(\bfv\otimes\bfzeta,\bfw\otimes\bfxi\big).
\end{align*}
\end{theorem}
\begin{proof}
We first remark that existence of the covariances on the right-hand sides follows from the Cauchy-Schwartz inequality in $L_{\mathbb{P}}^2$, \cref{lem:integrability}, and the assumptions on $\mu_1$ and $\varepsilon_1$.

To show the first equation we note that for all $\bfv,\bfw \in V^m$, $\bfzeta,\bfxi\in \mathbb{R}^m$, it follows from \cref{th:Saddlepoint} and \cref{lem:exp_bil} that
\begin{align*}
	\Cov&\begin{bmatrix}
	-a(\bfu_0,\cdot\,;\mu_1) \\ 0
	\end{bmatrix}\big(\bfv\otimes\bfzeta,\bfw\otimes\bfxi\big)\\  %
	&=
	\Cor
	\begin{bmatrix}
	A(D_{(0,0)}^{(1,0)} \bfu,\cdot)-B(D_{(0,0)}^{(1,0)} \bflambda,\cdot) \\ B(\cdot,D_{(0,0)}^{(1,0)} \bfu)
	\end{bmatrix}\big(\bfv\otimes\bfzeta,\bfw\otimes\bfxi\big)
	\\
	&= \int \limits_{\Omega} \begin{bmatrix}
	A(\cdot,\bfv)\otimes\Id & -B(\cdot,\bfv)\otimes\Id \\
	B(\bfzeta,\cdot)\otimes\Id & 0
	\end{bmatrix} \\ & \hspace*{2cm}
	\begin{bmatrix}
	D_{(0,0)}^{(1,0)} \bfu \otimes D_{(0,0)}^{(1,0)} \bfu 
	& D_{(0,0)}^{(1,0)} \bfu \otimes D_{(0,0)}^{(1,0)} \bflambda \\
	D_{(0,0)}^{(1,0)} \bflambda \otimes D_{(0,0)}^{(1,0)} \bfu 
	& D_{(0,0)}^{(1,0)} \bflambda \otimes D_{(0,0)}^{(1,0)} \bflambda
	\end{bmatrix} %
	\begin{bmatrix}
	\Id\otimes A(\cdot,\bfw) & \Id\otimes B(\bfxi,\cdot) \\
	-\Id\otimes B(\cdot,\bfw) & 0
	\end{bmatrix} \dif \mathbb{P} \\
	&= \begin{bmatrix} %
	A(\cdot,\bfv)\otimes\Id & -B(\cdot,\bfv)\otimes\Id \\
	B(\bfzeta,\cdot)\otimes\Id & 0
	\end{bmatrix}
	\Cor
	\begin{bmatrix}
	D_{(0,0)}^{(1,0)} \bfu  \\
	D_{(0,0)}^{(1,0)} \bflambda
	\end{bmatrix}
	\begin{bmatrix}
	\Id\otimes A(\cdot,\bfw) & \Id\otimes B(\bfxi,\cdot) \\
	-\Id\otimes B(\cdot,\bfw) & 0
	\end{bmatrix} 
\end{align*}
Finally, we note that
\[
\mathbb{E}[D_{(0,0)}^{(1,0)} \bfu]=0,
\qquad
\mathbb{E}[D_{(0,0)}^{(1,0)} \bflambda]=0,
\]
due to the same reasons as in \cref{th:stochastic}. Thus, the correlation in the last line is also a covariance and the assertion follows. The second equation follows in complete analogy.
\end{proof}
The extension of the theorem to uncentered perturbations is straightforward. We also note that \cref{eq:degen_der_eigs} also allows us to characterize the covariances of the eigenvalue derivatives directly. We state the following lemma without proof.
\begin{lemma}\label{lem:eigcovchar}
	Let $\lambda_0$ be an eigenvalue of multiplicity $m$ at $(\mu_0,\varepsilon_0)$ of \cref{eq:vareigprob} with $b(\cdot,\cdot\,;\varepsilon_0)$-orthonormal eigenbasis $\bfu_0$. Let $(\mu,\varepsilon)\mapsto(\bfu,\bflambda)$ be the unique local, analytic trajectory such that $(\bfu,\bflambda)$ satisfies \cref{eq:eig_bf,eq:eig_bf_normalized} with coefficients \cref{eq:perturbation_linear} and, at $(\mu_0,\varepsilon_0)$, it holds $(\bfu,\bflambda)=(\bfu_0,\lambda_0\bfI)$. Then it holds
\[
\Cov\Big[D_{(0,0)}^{(1,0)} \bflambda\Big]
=
\Cov\Big[a\big(\bfu_0,\bfu_0;\mu_1\big)\Big],\qquad
\Cov\Big[D_{(0,0)}^{(0,1)} \bflambda\Big]
=
\lambda_0^2\,\Cov\Big[b\big(\bfu_0,\bfu_0;\varepsilon_1\big)\Big].
\]
\end{lemma}
\subsection{On Karhunen-Lo\`eve-type random fields} \label{sec:KLE}
This subsection is concerned with the case that the stochastic variation $\mu_1$, $\varepsilon_1$ exhibit a Karhunen-Lo\`eve-type expansion, i.e., 
\begin{align} \label{eq:KLE}
	\mu_1(\omega) 
	= \sum \limits_{i=0}^{\infty} 
	z_{\mu,i}(\omega)\, L_{\mu,i} , \quad
	\varepsilon_1(\omega) 
	= \sum_{i=0}^{\infty} 
	z_{\varepsilon,i}(\omega)\, L_{\varepsilon,i}
\end{align}
with, for simplicity, $z_{\mu,i}, z_{\varepsilon,i}\sim \mathcal{U}[-\frac{1}{2},\frac{1}{2}]$ independent and uniformly distributed random variables and $L_{\mu,i} \in X$, $L_{\varepsilon,i} \in Y$. An immediate consequence is that the symmetry in of the random variables improves the accuracy of the covariance expansions from \cref{th:stochastic}.

\begin{corollary}\label{cor:symmetryexpansion}
For random fields of the form \cref{eq:KLE}, the covariance expansions from \cref{th:stochastic} are fourth order accurate, i.e., it holds
\begin{align*}		
\Cov[\bflambda] 
&= \alpha^2 \Cov[D_{(0,0)}^{(1,0)} \bflambda] 
+ \beta^2 \Cov[D_{(0,0)}^{(0,1)} \bflambda]
+ \LandauO\big((|\alpha|+|\beta|)^4\big), \\
\Cov[\bfu]
&= \alpha^2 \Cov[D_{(0,0)}^{(1,0)} \bfu] 
+ \beta^2 \Cov[D_{(0,0)}^{(0,1)} \bfu] 
+ \LandauO\big((|\alpha|+|\beta|)^4\big).
\end{align*}
\end{corollary}
\begin{proof}
In complete analogy to \cite[Lemma 2.3, 2.4, and 2.5]{HP2018a}.
\end{proof}

A further consequence is that the first order correction terms of the eigenvalue expansions can be characterized as follows.

\begin{lemma}
	Given the representation \cref{eq:KLE} for $\mu_1$ and $\varepsilon_1$ it holds
	\begin{align*}
	\Cov[D_{(0,0)}^{(1,0)} \bflambda] 
	&=\frac{1}{12}\sum_{k=1}^{\infty} 
	a(\bfu_0,\bfu_0;L_{\mu,k})\otimes a(\bfu_0,\bfu_0;L_{\mu,k}),\\  
	\Cov[D_{(0,0)}^{(0,1)} \bflambda] 
	&= \frac{\lambda^2_0}{12}\sum_{k=1}^{\infty} 
	 b(\bfu_0,\bfu_0;L_{\varepsilon,k}) \otimes b(\bfu_0,\bfu_0;L_{\varepsilon,k}),
	\end{align*}
in the statement of \cref{lem:eigcovchar}.
\end{lemma}
\begin{proof}
Follows from \cref{lem:eigcovchar} and the linearity of the bilinar forms in the parametric argument.
\end{proof}
The advantage of the lemma becomes obvious in implementations, when the series is truncated and the terms in the sum can be computed in parallel.

\section{Discretization} \label{sec:implementation}
\subsection{Galerkin discretization}
We briefly discuss the Galerkin discretization of the eigenvalue problem. For simplicity we restrict ourselves to the parameter model from \cref{sec:stochprob}. To that end, we assume a finite dimensional subspace $V_n\subset V$ spanned by basis functions $\{\varphi_1,\ldots,\varphi_n\}$ to be given and introduce the symmetric positive definite matrices
\begin{align*}
	\ul\bfA_\mu=\big[a(\varphi_i,\varphi_j;\mu)\big]_{i,j=1}^n, \qquad
	\ul\bfM_\varepsilon=\big[b(\varphi_i,\varphi_j;\varepsilon)\big]_{i,j=1}^n.
\end{align*}
This yields the discrete generalized eigenvalue problem for $1\leq m\leq n$ eigenpairs
\begin{align*}
	\text{find } (\ul\bfu ,\ul\bflambda) \in \mathbb{R}^{n\times m}\times \mathbb{R}^{m\times m} \text{ such that } 
	\ul\bfA_\mu \ul\bfu = \ul\bfM_\varepsilon \ul\bfu\, \ul\bflambda \text{ with } \ul\bflambda~\text{diagonal},
\end{align*}
to obtain approximate eigenvalues and coefficient vectors to approximate eigenfunctions to \cref{eq:vareigprob}. The precise approximation properties depend on the approximation properties of the subspace $V_n$, see, e.g., \cite{Bof2010}.

\subsection{Eigenpair derivatives}\label{sec:disceigpairder}
Having computed a Galerkin approximation $(\ul\bfu_0,\ul\bflambda_0)$ to the reference solution $(\bfu_0,\bflambda_0)$ of \cref{eq:vareigprob}, the saddle point problem for the derivatives \cref{eq:saddlepoint_bfu} can also be computed by the Galerkin method. Using  Galerkin discretization of \cref{eq:saddlepoint_bfu} yields the discrete saddle point problems
\begin{align*}
\begin{bmatrix}
\ul\bfA_{\mu_0}-\ul\bflambda_0 \,\ul\bfM_{\varepsilon_0} & -\ul\bfM_{\varepsilon_0} \ul\bfu_0 \\
\ul\bfu_0^\intercal \ul\bfM_{\varepsilon_0} & \bf0
\end{bmatrix}
\begin{bmatrix}
\ul{D_{(0,0)}^{(1,0)} \bfu} \\
\ul{D_{(0,0)}^{(1,0)} \bflambda}
\end{bmatrix}
&= 
\begin{bmatrix}
-\ul\bfA_{\mu_1} \ul\bfu_0 \\
\bf0 
\end{bmatrix} ,
\\ %
\begin{bmatrix}
\ul\bfA_{\mu_0}-\ul\bflambda_0 \,\ul\bfM_{\varepsilon_0} & -\ul\bfM_{\varepsilon_0} \ul\bfu_0 \\
\ul\bfu_0^\intercal \ul\bfM_{\varepsilon_0} & \bf0
\end{bmatrix}
\begin{bmatrix}
\ul{D_{(0,0)}^{(0,1)} \bfu} \\
\ul{D_{(0,0)}^{(0,1)} \bflambda}
\end{bmatrix}
&= 
\begin{bmatrix}
\ul\bfM_{\varepsilon_1} \ul\bfu_0 \ul\bflambda_0 \\ 
-\frac{1}{2}\diag_{i=1,\ldots,m}([\ul\bfu_{0}]_i^\intercal \ul{\bfM}_{\varepsilon_1} [\ul\bfu_{0}]_i)
\end{bmatrix}.
\end{align*}
However, a few remarks are in order. First, unique solvability of the systems can be shown in complete analogy to the continuous case by a discrete LBB-condition. Second, the Galerkin approximations to the analytical eigenpairs are required for assembly of the system. This leads to a consistency error whose analysis is outside the scope of this paper, but is likely to be manageable when considering the precise approximation properties of $V_n$ combined with Strang's lemma. Third, the obtained discrete system is the same as if we would apply \cref{eq:saddlepoint_bfu} to the discrete system. Thus, the ``derive and then discretize'' and the ``discretize and then derive'' approach can be considered equivalent up to the addressed consistency error.

Of course, in analogy to \cref{eq:degen_der_eigs}, the approximate eigenvalue derivatives can also be characterized due to
\[
\ul{D_{(0,0)}^{(1,0)} \bflambda}
=
\ul\bfu_0^\intercal\ul\bfA_{\mu_1}\ul\bfu_0,
\qquad
\ul{D_{(0,0)}^{(0,1)} \bflambda}
=
-\ul\bfu_0^\intercal\ul\bfM_{\varepsilon_1}\ul\bfu_0\ul\bflambda_0.
\]

\subsection{Covariance equations}
Discretizing the characterization of covariances of eigenpair derivatives from \cref{th:TensorVarDer} by means of the Galerkin method yields the matrix equations
\begin{subequations}\label{eq:disccoveqs}
\begin{align}
\begin{bmatrix}
\ul\bfA_{\mu_0}-\ul\bflambda_0 \, \ul\bfM_{\varepsilon_0} & -\ul\bfM_{\varepsilon_0} \ul\bfu_0 \\
\ul\bfu_0^\intercal \ul\bfM_{\varepsilon_0} & \bf0
\end{bmatrix}
\Cov
\begin{bmatrix}
\ul{D_{(0,0)}^{(1,0)} \bfu} \\
\ul{D_{(0,0)}^{(1,0)} \bflambda}
\end{bmatrix}
\begin{bmatrix}
\ul\bfA_{\mu_0}-\ul\bflambda_0 \, \ul\bfM_{\varepsilon_0} & \ul\bfM_{\varepsilon_0} \ul\bfu_0 \\
-\ul\bfu_0^\intercal \ul\bfM_{\varepsilon_0} & \bf0
\end{bmatrix}
&=
\Cov
\begin{bmatrix}
-\ul\bfA_{\mu_1} \ul\bfu_0 \\
\bf0
\end{bmatrix},\\
\begin{bmatrix}
\ul\bfA_{\mu_0}-\ul\bflambda_0 \, \ul\bfM_{\varepsilon_0} & -\ul\bfM_{\varepsilon_0} \ul\bfu_0 \\
\ul\bfu_0^\intercal \ul\bfM_{\varepsilon_0} & \bf0
\end{bmatrix}
\Cov
\begin{bmatrix}
\ul{D_{(0,0)}^{(0,1)} \bfu} \\
\ul{D_{(0,0)}^{(0,1)} \bflambda}
\end{bmatrix}
\begin{bmatrix}
\ul\bfA_{\mu_0}-\ul\bflambda_0 \, \ul\bfM_{\varepsilon_0} & \ul\bfM_{\varepsilon_0} \ul\bfu_0 \\
-\ul\bfu_0^\intercal \ul\bfM_{\varepsilon_0} & \bf0
\end{bmatrix}\\
&\hspace*{-2.1cm}=
\Cov
\begin{bmatrix}
\ul\bfM_{\varepsilon_1} \ul\bfu_0\, \ul\bflambda_0 \\ 
-\frac{1}{2}\diag_{i=1,\ldots,m}([\ul\bfu_{0}]_i^\intercal \ul{\bfM}_{\varepsilon_1} \ul [\bfu_{0}]_i)
\end{bmatrix},
\end{align}
\end{subequations}
for which unique solvability follows from the considerations in \cref{sec:disceigpairder}. Unfortunately, covariance matrices are usually densely populated such that naive approaches to solve the above equations are prohibitively expensive to solve for approximation spaces of sufficiently many degrees of freedom, even if $\ul\bfA_{\mu_0}$ and $\ul\bfB_{\varepsilon_0}$ are sparse. Fortunately, quite a few articles have addressed the efficient solution of such systems by various means such as sparse grids, global low-rank approximations, hierarchical matrices, and others \cite{BG2004a,DHP2017,DHS2017,HPS2012,HSS2008,vS2006}. In our numerical experiments below, we use a global low-rank approach where the covariance matrices on the right-hand side are approximated by low-rank matrices. Given access to on-the-fly computable matrix entries, the algorithm provides a black-box strategy to obtain an error-controlled low-rank approximation of rank $k\ll n$ without the full assembly of the covariance matrices. A low-rank approximation to the solution covariances can then straightforwardly obtained by solving $k$ systems of linear equations, possibly in parallel. Of course, if a (finite dimensional) Karhunen-Lo\`eve type expansion of the random fields such as in \cref{eq:KLE} is available, we may use this decomposition directly.

\section{Numerical Examples} \label{sec:experiments}

\subsection{Problem setting}
The numerical examples to illustrate our findings are based on the diffusion equation from \cref{sec:examples} on the unit square $D=(0,1)^2$. 
To that end, the problem is discretized by continuous, piecewise linear finite elements with $n_{\text{FE}} = 481$ degrees of freedom. 
Here, we note that the limiting factor are the computational cost of the Monte Carlo method, which serves as a reference solution. 
The perturbation approach can easily deal with more degrees of freedom. 

The coefficients $\mu$ and $\varepsilon$ are modelled as truncated Karhunen-Lo\`eve expansions, i.e., as in \cref{eq:perturbation_linear} with $\mu_0=\varepsilon_0=1$ and $\mu_1$ and $\varepsilon_1$ as in \cref{eq:KLE} with the series truncated after finitely many terms and $L_{\mu,i}=L_{\varepsilon,i}=\sqrt{\sigma_i}\phi_i$. Here, $(\phi_i,\sigma_i)$ corresponds to the $i$-th eigenpair of the covariance operator
\[
\mathcal{C}\colon L^2(D)\to L^2(D),
\qquad
(\mathcal{C}\phi)(\bfx)=\int_Dg(\bfx,\bfy)\phi(\bfy)\dd\bfy,
\quad
\bfx\in D,
\]
with covariance kernel
\begin{align*}
g(\bfx,\bfy)=g(\|\bfx-\bfy\|_2),
\qquad
g(r) = \frac{1}{\sqrt{20\pi}}\exp\Big(-\frac{r^2}{20}\Big).
\end{align*}
To this end, to actually compute the eigenpairs of the covariance operator numerically, we discretize the eigenvalue problem $\mathcal{C}\phi_i=\sigma_i\phi_i$ by means of the Galerkin method using continuous and piecewise linear finite elements to obtain a discrete eigenvalue problem $\underline{\bfC}\underline{\bfphi}_i=\underline{\sigma}_i\bfM\underline{\bfphi}_i$. The discrete eigenvalue problem is then solved approximately by means of the pivoted Cholesky decomposition as follows, see also \cite{HPS2012, HPS2015} for more details and error estimates. After adaptively computing a rank-$k$ factorization $\underline{\bfC}\approx\underline{\widetilde{\bfL}}_k\underline{\widetilde{\bfL}}_k^\intercal$ up to a tolerance of $10^{-5}$ in the trace norm by means of the pivoted Cholesky decomposition, we solve the $k\times k$ eigenvalue problem
\[
\underline{\widetilde{\bfL}}_k^\intercal\underline{\bfM}^{-1}\underline{\widetilde{\bfL}}_k\underline{\widetilde{\bfphi}}_i=\widetilde{\sigma}_i\underline{\widetilde{\bfphi}_i}.
\]
In our case, the required rank determined by the pivoted Cholesky decomposition was $k=283$, such that the solution of the eigenvalue problem can be accomplished by standard dense linear algebra. Finally, $(\underline{\bfphi}_i,\underline{\sigma}_i)\approx(\underline{\bfM}^{-1}\underline{\bfL}\widetilde{\bfphi}_i,\widetilde{\sigma}_i)$ is a sufficient approximation to the $283$ eigenpairs with the largest eigenvalues of $\mathcal{C}$. This yields approximations of the form \cref{eq:KLE} which are truncated after $283$ terms.

\begin{figure}[]
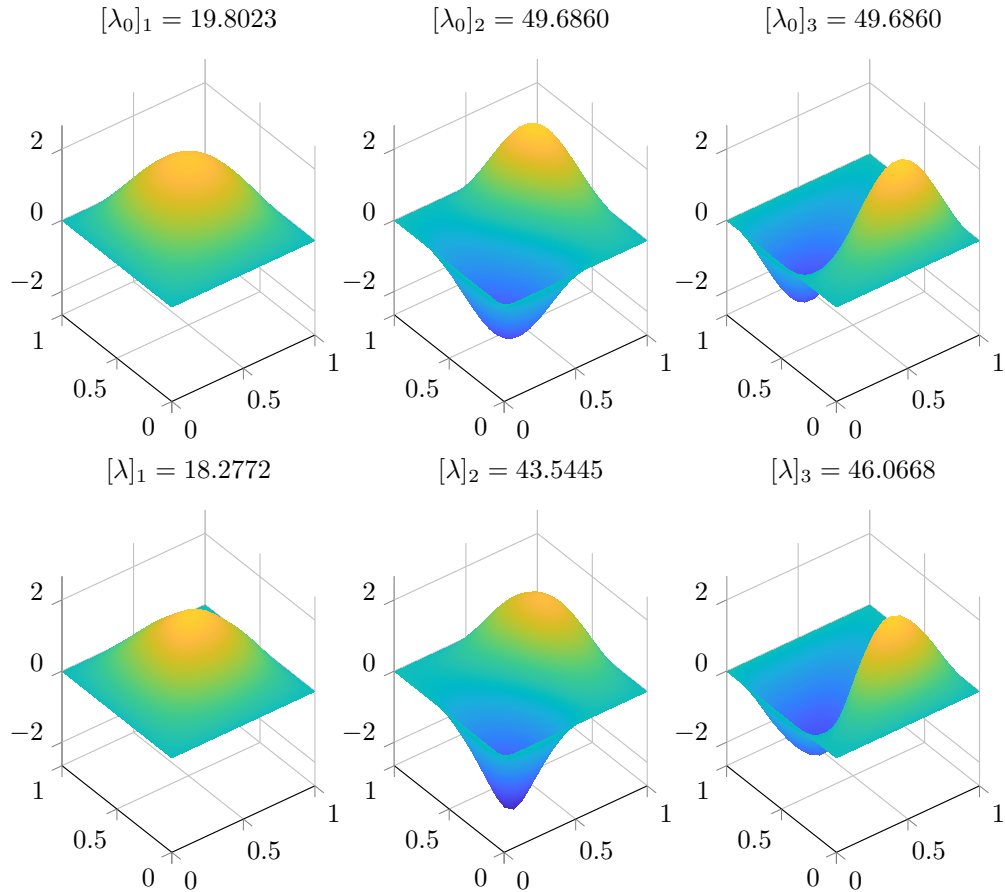

	\center
	\input{./fig/fig_eigenfunctions.tex}
	\input{./fig/fig_eigenfunctions_pert.tex}
	\caption{
	Unperturbed (top) and perturbed (bottom) first three eigenpairs of the diffusion operator for $\alpha=\beta=1$. The corresponding eigenvalue trajectories are illustrated in \cref{fig:splitting_eigvalue}.
	} 
	\label{fig:eigenfunctions}
\end{figure}

\begin{figure}[]
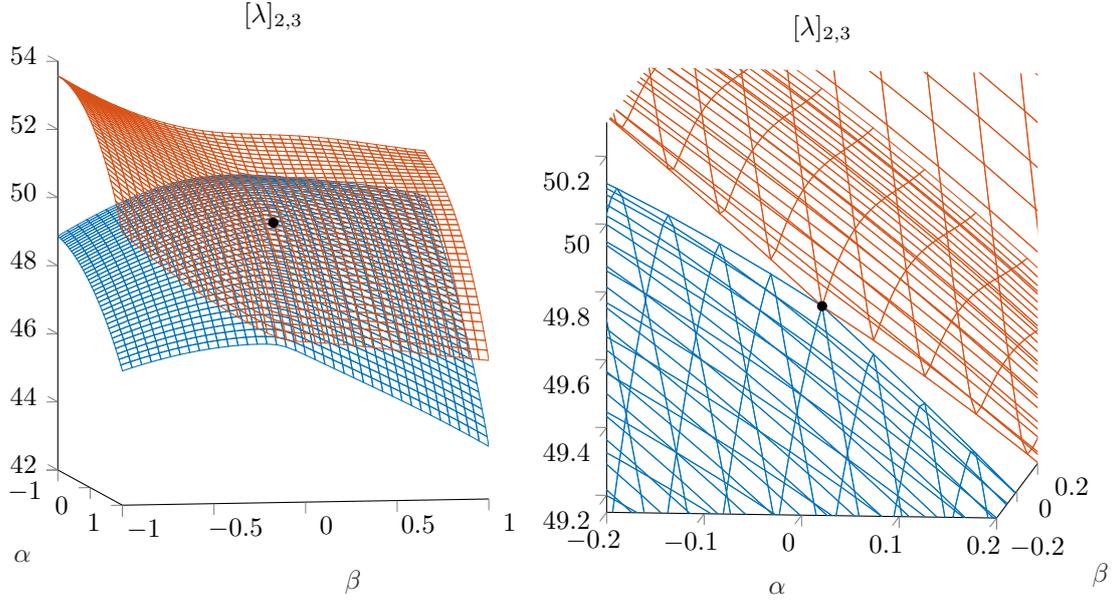

	\center
	\input{./fig/fig_splitting_eigvalue.tex}
	\input{./fig/fig_splitting_eigvalue_2.tex}
	\caption{
	Illustration of the eigenvalue trajectories of $[\lambda]_{2,3}$ under a sample perturbation.
	} \label{fig:splitting_eigvalue}
\end{figure}

It is well known that the diffusion equation on the unit square has a unique first eigenvalue, whereas the second and third eigenvalue coincide, yielding a degenerate eigenvalue of multiplicity $m=2$. We will consider the first eigenspace as a nondegenerate and the second eigenspace as a degenerate eigenspace both in isolation. An illustration of the three eigenpairs for the unperturbed case and for a sample perturbation can be found in \cref{fig:eigenfunctions}.
Looking at the corresponding eigenvalues in the titles, it becomes clear that the multiple eigenvalue of the unperturbed problem splits with the perturbation. The precise trajectories of the eigenvalues for this sample perturbation are illustrated in \cref{fig:splitting_eigvalue}.

\subsection{Deterministic approximation of eigenpair through expansions}\label{sec:num_det}
We confirm the eigenpair approximations due to the deterministic expansions from \cref{sec:deterministic}. To this end, we fix a realization of $\mu_1$ and $\varepsilon_1$ for all $\alpha,\beta\in \{2^i,i\in\{-15,-14,\ldots,0\}\}$.
From the illustration in \cref{fig:splitting_eigvalue} it becomes clear that the trajectories form a cone-like shape, such that a common polarization matrix can never be found for this sample perturbation, see also \cref{rem:rewritefrechet}.
To align the reference eigenpair and the perturbed eigenpair for calculation of the error comparison we compare two approaches. The first approach is to compute a polarization for each $(\alpha,\beta)$ (cf. \cref{sec:polarization}). The second approach is to rotate the eigenspaces onto each other by means of an SVD (cf. \cref{sec:sampling}).

\begin{figure}[]
	\center
	\input{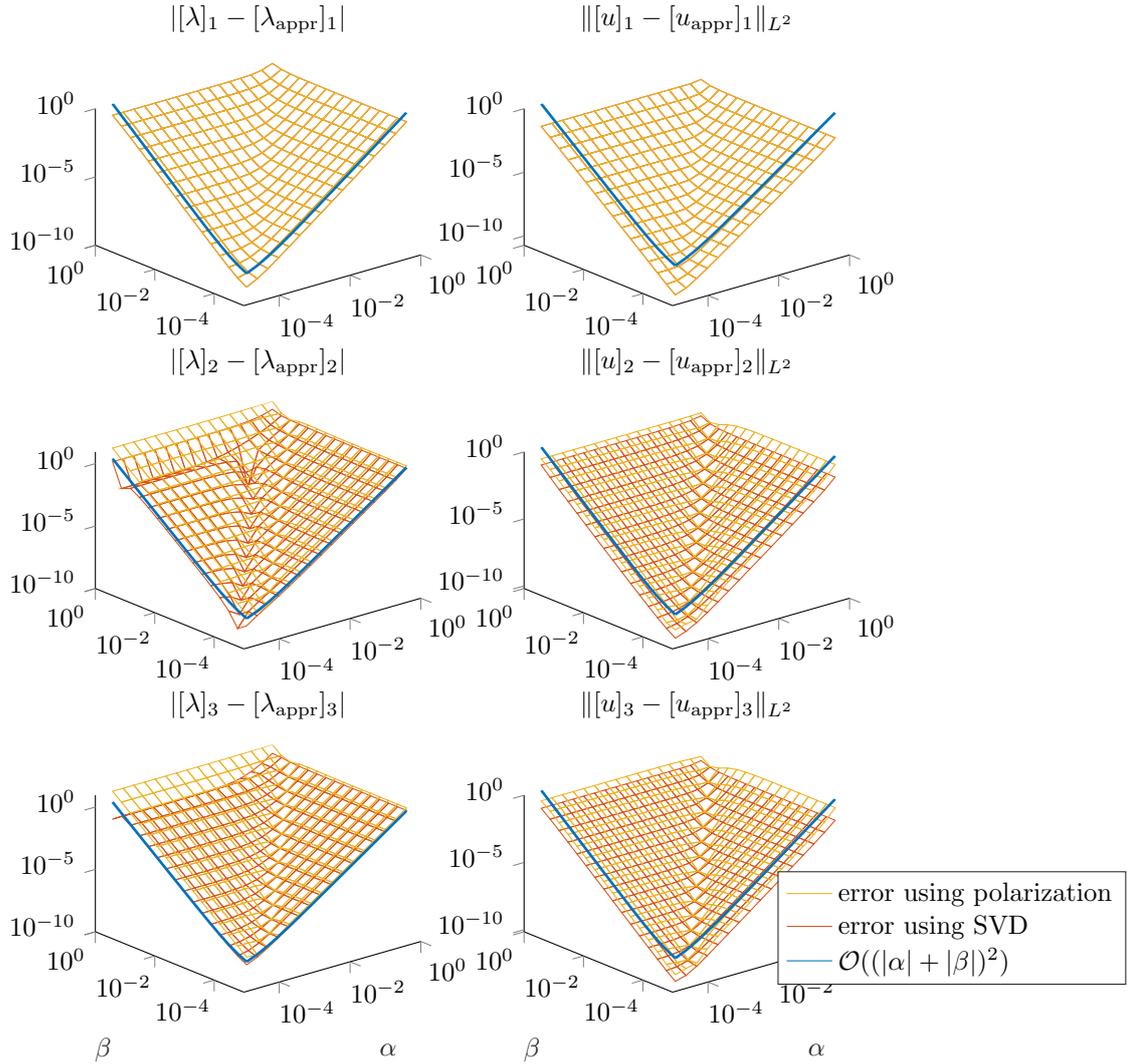}
	\caption{Convergence of residue of first-order approximations in the determinstic setting.
	} \label{fig:determ}
\end{figure}

From the theoretical considerations in \cref{sec:deterministic} we expect an approximation error of $\LandauO((|\alpha|+|\beta|)^2)$, which can clearly be confirmed from the error graphs in \cref{fig:determ}.

\subsection{Uncertainty quantification using Monte Carlo} \label{sec:UQ_MC}
As outlined in \cref{sec:behav_stoch_deg}, polarizations are not necessarily meaningful in the context of stochastic perturbations and looking at the eigenspaces seems to be more reasonable. In the following, we confirm convergence of the Monte Carlo method for mean and covariance when the eigenspaces are aligned as described in \cref{sec:sampling}.

To this end, to estimate the root mean square error of the Monte Carlo method,
we exploit that for a $W$-valued random variable $w\in L_{\mathbb{P}}^2(\Omega;W)$ and its Monte Carlo estimator with $M$ samples $E_M[w]$ it holds
\[
\|\mathbb{E}[w]-E_M[w]\|_{L_{\mathbb{P}}^2(\Omega;W)}^2
=
\frac{1}{M}\mathbb{E}\big[\|\mathbb{E}[w]-w\|_W^2\big]
\approx
\frac{1}{M^2}\sum_{i=1}^M\bigg\|\bigg(\frac{1}{M}\sum_{j=1}^Mw_j\bigg)-w_i\bigg\|_W^2.
\]
For $W=\mathbb{R}$, the latter quantity is the sample variance scaled with $1/M$ of the samples from the Monte Carlo estimator, which is frequently used to estimate errors from Monte Carlo estimators. In our case, it holds $W=\mathbb{R}^{m\times m}$ for eigenvalues and $W=L^2(D)$ for eigenfunctions for the mean of eigenpairs with multiplicity $m$. For the mean square error of the covariance, the same considerations hold with $W=\mathbb{R}^{m\times m}\otimes \mathbb{R}^{m\times m}$ for eigenvalues and $W=L^2(D\times D)$ for eigenfunctions.

For the Monte Carlo estimator of the mean standard sampling as well as antithetic sampling was tested. This is justified, since the approximately linear nature of the expected value (cf. \cref{th:stochastic}) leads to a lower expected error when using antitethic sampling.
In order to accurately calculate the error for the antithetic sampling, the average over the antithetic pair was counted as two samples. 
Since the antithetic scheme was highly effective for the chosen perturbation, the antithetic estimate of the expected values was also used for centering the sample of the covariances, which were otherwise chosen according to standard Monte Carlo, since the covariance is approximately quadratic for small perturbations. 

\begin{figure}
	\center
	\definecolor{mycolor1}{rgb}{0.00000,0.44706,0.74118}%
\definecolor{mycolor2}{rgb}{0.85098,0.32549,0.09804}%
\begin{tikzpicture}

\begin{axis}[%
width=0.38\textwidth,
height=0.4\textwidth,
at={(0\textwidth,0\textwidth)},
scale only axis,
xmode=log,
xmin=100,
xmax=10000000,
xminorticks=true,
xlabel style={font=\color{white!15!black}},
xlabel={$M$, number of samples},
ymode=log,
ymin=1e-16,
ymax=1e-1,
yminorticks=true,
axis background/.style={fill=white},
title style={font=\bfseries},
title={$\bf\lambda$, $\|\cdot\|_{L_{\mathbb{P}}^2(\Omega;W)}^2$-error},
axis x line*=bottom,
axis y line*=left,
legend style={legend cell align=left, align=left, draw=white!15!black, font=\tiny, at={(0,0)},anchor=south west},
legend columns=2,
domain=1e2:1e7
]
\addplot [color=mycolor1, dotted, mark=x, mark options={solid, mycolor1}]
  table[row sep=crcr]{%
100	4.3772375743925e-08\\
1000	2.94186345545942e-09\\
10000	2.6825409603056e-10\\
100000	2.72115891126432e-11\\
1000000	2.72488430948059e-12\\
4000000	6.78997212524959e-13\\
8000000	3.39634812180855e-13\\
};
\addlegendentry{$\mathbb{E}\big[[\bflambda]_1\big]$, antit.}

\addplot [color=mycolor1, dotted, mark=+, mark options={solid, mycolor1}]
  table[row sep=crcr]{%
100	0.000225996514947292\\
1000	2.0250915158964e-05\\
10000	2.01106642013106e-06\\
100000	2.01353086801255e-07\\
1000000	2.02187415332813e-08\\
4000000	5.05578592330168e-09\\
8000000	2.52643083963506e-09\\
};
\addlegendentry{$\mathbb{E}\big[[\bflambda]_1\big]$, std.}

\addplot [color=mycolor1, dotted, mark=asterisk, mark options={solid, mycolor1}]
  table[row sep=crcr]{%
100	1.35170536752956e-05\\
1000	9.73442200565423e-07\\
10000	9.99169247207378e-08\\
100000	9.9554037077856e-09\\
1000000	1.00132882539751e-09\\
4000000	2.50423108042741e-10\\
8000000	1.25033632940512e-10\\
};
\addlegendentry{$\operatorname{Cov}\big[[\bflambda]_1\big]$, std.}

\addplot [color=mycolor2, dotted, mark=x, mark options={solid, mycolor2}]
  table[row sep=crcr]{%
100	7.87583836625575e-07\\
1000	7.61630660245272e-08\\
10000	7.77437892830145e-09\\
100000	7.7350361062383e-10\\
1000000	7.74159343500115e-11\\
4000000	1.93344453488171e-11\\
8000000	9.66117309334801e-12\\
};
\addlegendentry{$\mathbb{E}\big[[\bflambda]_{2,3}\big]$, antit.}

\addplot [color=mycolor2, dotted, mark=+, mark options={solid, mycolor2}]
  table[row sep=crcr]{%
100	0.00250358836286263\\
1000	0.000228956055293018\\
10000	2.27348456252609e-05\\
100000	2.282873937275e-06\\
1000000	2.29186297451021e-07\\
4000000	5.73138262361931e-08\\
8000000	2.86494052112032e-08\\
};
\addlegendentry{$\mathbb{E}\big[[\bflambda]_{2,3}\big]$, std.}

\addplot [color=mycolor2, dotted, mark=asterisk, mark options={solid, mycolor2}]
  table[row sep=crcr]{%
100	0.00109990155898258\\
1000	9.41336933275294e-05\\
10000	9.19926836961741e-06\\
100000	9.19253477747575e-07\\
1000000	9.26214008065463e-08\\
4000000	2.31822717490838e-08\\
8000000	1.15792280740925e-08\\
};
\addlegendentry{$\operatorname{Cov}\big[[\bflambda]_{2,3}\big]$, std.}

\addplot[dashed] {1/x};
\addlegendentry{$\mathcal{O}(1/M)$}

\end{axis}
\end{tikzpicture}%
	\definecolor{mycolor1}{rgb}{0.00000,0.44706,0.74118}%
\definecolor{mycolor2}{rgb}{0.85098,0.32549,0.09804}%
\begin{tikzpicture}

\begin{axis}[%
width=0.38\textwidth,
height=0.4\textwidth,
at={(0\textwidth,0\textwidth)},
scale only axis,
xmode=log,
xmin=100,
xmax=10000000,
xminorticks=true,
xlabel style={font=\color{white!15!black}},
xlabel={$M$, number of samples},
ymode=log,
ymin=1e-16,
ymax=1e-1,
yminorticks=true,
axis background/.style={fill=white},
title style={font=\bfseries},
title={$\bfu$, $\|\cdot\|_{L_{\mathbb{P}}^2(\Omega;W)}^2$-error},
axis x line*=bottom,
axis y line*=left,
legend style={legend cell align=left, align=left, draw=white!15!black, font=\tiny, at={(0,0)},anchor=south west},
legend columns=2,
domain=1e2:1e7
]
\addplot [color=mycolor1, dotted, mark=x, mark options={solid, mycolor1}]
  table[row sep=crcr]{%
100	8.99362791234205e-09\\
1000	9.70021541861737e-10\\
10000	1.08756258967517e-10\\
100000	1.09695719251036e-11\\
1000000	1.09387151335533e-12\\
4000000	2.72942042750051e-13\\
8000000	1.36442230709963e-13\\
};
\addlegendentry{$\mathbb{E}\big[[\bfu]_1\big]$, antit.}

\addplot [color=mycolor1, dotted, mark=+, mark options={solid, mycolor1}]
  table[row sep=crcr]{%
100	2.82643175277931e-05\\
1000	2.83720711162895e-06\\
10000	2.76913596237088e-07\\
100000	2.76730355186119e-08\\
1000000	2.76484756428497e-09\\
4000000	6.91156373719914e-10\\
8000000	3.45551639195074e-10\\
};
\addlegendentry{$\mathbb{E}\big[[\bfu]_1\big]$, std.}

\addplot [color=mycolor1, dotted, mark=asterisk, mark options={solid, mycolor1}]
  table[row sep=crcr]{%
100	1.02454019052217e-07\\
1000	9.98944822627621e-09\\
10000	9.45479865543844e-10\\
100000	9.45617375148832e-11\\
1000000	9.43670308169546e-12\\
4000000	2.35883695608289e-12\\
8000000	1.17906954430606e-12\\
};
\addlegendentry{$\operatorname{Cov}\big[[\bfu]_1\big]$, std.}

\addplot [color=mycolor2, dotted, mark=x, mark options={solid, mycolor2}]
  table[row sep=crcr]{%
100	1.82880573315292e-07\\
1000	1.96001574676002e-08\\
10000	2.10248378560601e-09\\
100000	2.10153195704862e-10\\
1000000	2.09359255951161e-11\\
4000000	5.23087500772771e-12\\
8000000	2.6139977162243e-12\\
};
\addlegendentry{$\mathbb{E}\big[[\bfu]_{2,3}\big]$, antit.}

\addplot [color=mycolor2, dotted, mark=+, mark options={solid, mycolor2}]
  table[row sep=crcr]{%
100	0.000288407858869928\\
1000	2.86360488784997e-05\\
10000	2.83436911223782e-06\\
100000	2.84101571218024e-07\\
1000000	2.84144620730812e-08\\
4000000	7.10399598528752e-09\\
8000000	3.55208543069146e-09\\
};
\addlegendentry{$\mathbb{E}\big[[\bfu]_{2,3}\big]$, std.}

\addplot [color=mycolor2, dotted, mark=asterisk, mark options={solid, mycolor2}]
  table[row sep=crcr]{%
100	1.03819093546673e-05\\
1000	1.00727090859803e-06\\
10000	9.84910688397305e-08\\
100000	9.89630846742658e-09\\
1000000	9.89825328083452e-10\\
4000000	2.47489381376728e-10\\
8000000	1.23751904414397e-10\\
};
\addlegendentry{$\operatorname{Cov}\big[[\bfu]_{2,3}\big]$, std.}

\addplot[dashed] {1/x};
\addlegendentry{$\mathcal{O}(1/M)$}

\end{axis}
\end{tikzpicture}%
	\caption{
	Convergence of Monte Carlo method for eigenvalues and -spaces.
	} 
	\label{fig:MC}
\end{figure}
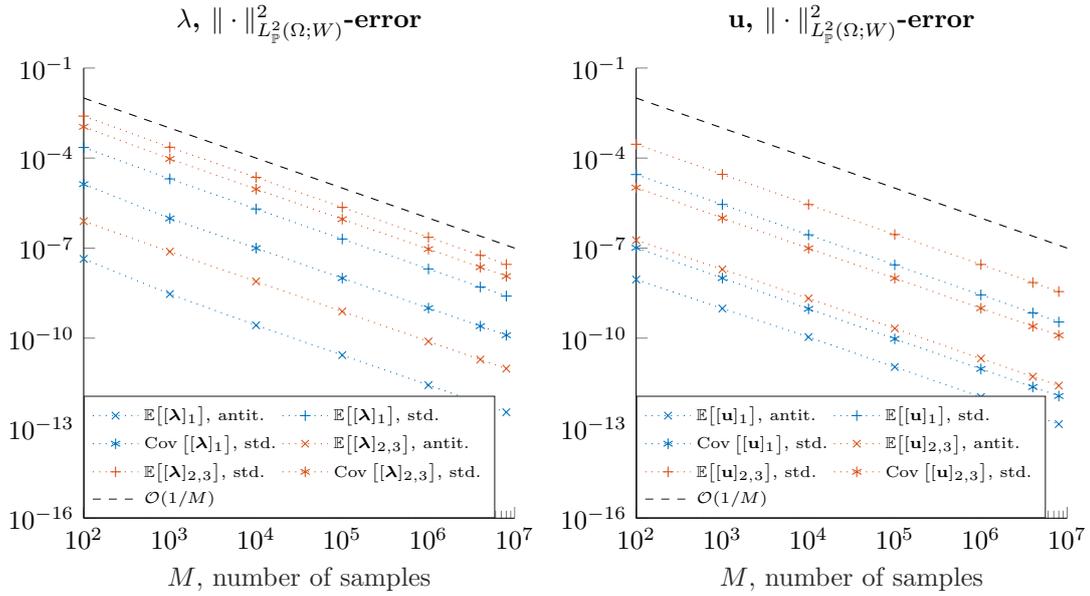

The graphs in \cref{fig:MC} confirm the expected order of convergence. To that end, we used the $L^2$-norm for the eigenfunctions and the Frobenius norm for the eigenvalues.

\subsection{Uncertainty quantification using expansions}
It remains to verify the expansions for mean and covariance from \cref{th:stochastic}. To that end, we take the same combinations for $\alpha$ and $\beta$ as for the examples for the deterministic expansions in \cref{sec:num_det} and compute for each $(\alpha,\beta)$ combination a Monte Carlo reference solution with $4 \cdot 10^6$ samples for mean and covariance. 
The covariance equations \cref{eq:disccoveqs} for the second order correction terms of $\Cov[\bfu]$ and $\Cov[\bflambda]$ are solved by means of the low-rank approach. To this end, low-rank approximations of the right-hand sides of the covariance equations \cref{eq:disccoveqs} are given through
\begin{align*}
\Cov
\begin{bmatrix}
-\ul\bfA_{\mu_1} \ul\bfu_0 \\
\bf0
\end{bmatrix}
&\approx
\frac{1}{12}\sum_{i=1}^k
\begin{bmatrix}
\ul{\bfA}_{L_{\mu,i}}\ul{\bfu}_0 \\
\bf0 
\end{bmatrix}
\begin{bmatrix}
\big(\ul{\bfA}_{L_{\mu,i}}\ul{\bfu}_0\big)^\intercal &
\bf0 
\end{bmatrix},\\
\Cov
\begin{bmatrix}
\ul\bfM_{\varepsilon_1} \ul\bfu_0\, \ul\bflambda_0 \\ 
\ul\bfN_{\varepsilon_1}
\end{bmatrix}
&\approx
\frac{1}{12}\sum_{i=1}^k
\begin{bmatrix}
\ul{\bfM}_{L_{\varepsilon,i}}\ul{\bfu}_0\ul\bflambda_0 \\
\ul\bfN_{L_{\varepsilon,i}}
\end{bmatrix}
\begin{bmatrix}
\big(\ul{\bfM}_{L_{\varepsilon,i}}\ul{\bfu}_0\ul\bflambda_0\big)^\intercal &
\ul\bfN_{L_{\varepsilon,i}}^\intercal
\end{bmatrix},
\end{align*}
where we abbreviate $\ul\bfN_{\varepsilon_1}=-\frac{1}{2}\diag_{i=1,\ldots,m}([\ul\bfu_{0}]_i^\intercal \ul{\bfM}_{\varepsilon_1} \ul [\bfu_{0}]_i)$. Substituting these low-rank factorizations into the covariance equations \cref{eq:disccoveqs}, a low-rank approximation of
\[
\Cov
\begin{bmatrix}
\ul{D_{(0,0)}^{(1,0)} \bfu} \\
\ul{D_{(0,0)}^{(1,0)} \bflambda}
\end{bmatrix},
\qquad
\Cov
\begin{bmatrix}
\ul{D_{(0,0)}^{(0,1)} \bfu} \\
\ul{D_{(0,0)}^{(0,1)} \bflambda}
\end{bmatrix},
\]
can be computed by solving $k$ linear systems of equations, see, e.g., \cite[Chapter 4.1]{HPS2012}.

The convergence results for the mean are shown in \cref{fig:exp} and in line with the results from \cref{th:stochastic}. The convergence results for the covariance are shown in \cref{fig:cov_res} and seem to be as expected from \cref{th:stochastic} and \cref{cor:symmetryexpansion}. 
To that end, we may observe for the covariances that the errors for $(\alpha,\beta)$ close to the origin are expected to be near machine precision, which is hard to achieve with a Monte Carlo reference solution. 
It seems reasonable to assume that the reduced order close to the origin is due to round-off errors and the missing accuracy of the reference solution.

\begin{figure}[p]
	\center
	\input{./fig/fig_exp.tex}
	\caption{
	Convergence of series approximations of the mean. 
    }
	\label{fig:exp}
	\center
	\input{./fig/fig_cov_res.tex}
	\caption{Convergence of series approximations of the covariance.
   } 
	\label{fig:cov_res}
\end{figure}

\subsection{A remark on computational cost}

All numerical examples are implemented in Matlab 2022a and run on a compute server with four Intel(R) Xeon(R) E7-4850 v2 CPU with twelve 2.30GHz cores each and hyperthreading disabled and 1.3 TB RAM. The Monte Carlo computations were parallelized over the number of samples and the solution of the covariance equations over the terms of the Karhunen-Lo\`eve expansion. All sparse matrices were stored as such and computations were done using sparse linear algebra, except when explicitly mentioned.

The complexity of computing a rank-$k$ approximation by means of the pivoted Cholesky decomposition is $\mathcal{O}(nk^2)$. Solving the $k\times k$ eigenvalue problem has a complexity of $\mathcal{O}(k^3)$, which allows to compute a Karhunen-Lo\`eve expansion in $\mathcal{O}(nk^2)$ operations. Likewise, assuming that the solution of a sparse saddle point problem can be accomplished in $\mathcal{O}(n)$ operations, the solution of the covariance equations can be accomplished in $\mathcal{O}(nk^2)$ operations.

To accelerate the assembly of the system matrices throughout the Monte Carlo simulation, the matrices $\underline{\bfA}_{L_{\mu,i}},\underline{\bfM}_{L_{\varepsilon,i}}$ are assembled a priori. This allows for an accelerated assembly of $\ul\bfA_\mu$ and $\ul\bfM_\varepsilon$ through linear combinations throughout the sampling process, otherwise the Monte-Carlo estimate would be additionally disadvantaged by repeated discretization of the random field sample. The matrices $\underline{\bfA}_{L_{\mu,i}},\underline{\bfM}_{L_{\varepsilon,i}}$ are also required for the perturbation approach. The solution time of the eigenvalue problem at the reference point $(\alpha,\beta)=(0,0)$ was negligible for the perturbation approach. The solution of the covariance equations took 14 seconds. Thus, the perturbation approach offers approximations for mean and covariance of the eigenpairs for all combinations of $(\alpha,\beta)$ after only a few seconds. In contrast, the sampling process alone of the Monte Carlo approach for our reference solution with $4\cdot 10^6$ samples required approximately 35 minutes for a single $(\alpha,\beta)$ pairing. Thus, the computations for all reference points in \cref{fig:exp} and \cref{fig:cov_res} took more than 6 days. It should be clear that comparisons with larger finite element spaces are beyond our reach at the moment.

\section{Conclusion} \label{sec:conclusion}
We considered uncertainty quantification approaches to generalized eigenvalue problems with stochastic parameter dependence. To that end, our main interest was on eigenpairs with higher but finite multiplicity, where crossings and bifurcations of the eigenvalue trajectories are possible. To improve our understanding of the situation we considered the Fr\'echet derivatives of the eigenpairs and provided a new \emph{linear} characterization for the derivatives and a new regularity result of the eigenpairs. With this improved understanding we concluded that the uncertainty quantification of eigenvalues with higher multiplicity is not meaningful in general and that the uncertainty quantification of eigenspaces seems to be a reasonable approach. We discussed a strategy how to relate eigenpairs of different samples for sampling-based methods and a perturbation approach for the uncertainty quantification of eigenpairs. We provided numerical examples for both approaches, illustrating their feasibility.

Finally, we would like to remark that the same strategy as for deriving the first derivatives in \cref{sec:deterministic} could also be used to derive higher Fr\`echet derivatives of the eigenpairs. These derivatives could be used to derive more accurate expansions for the mean, covariance and other statistical quantities of interest, see also \cite{DH2018,Doe2020} for a discussion.

\bibliographystyle{plain}
\bibliography{bibliography}
\end{document}